# The asymptotically optimal estimating equation for longitudinal data. Strong Consistency

R.M. Balan[*]  L. Dumitrescu[†]
University of Ottawa  University of Ottawa

I. Schiopu-Kratina [‡]
Statistics Canada

July 4, 2008

**Abstract**

In this article, we introduce a conditional marginal model for longitudinal data, in which the residuals form a martingale difference sequence. This model allows us to consider a rich class of estimating equations, which contains several estimating equations proposed in the literature. A particular sequence of estimating equations in this class contains a random matrix $\mathcal{R}_{i-1}^*(\beta)$, as a replacement for the "true" conditional correlation matrix of the $i$-th individual. Using the approach of [12], we identify some sufficient conditions under which this particular sequence of equations is asymptotically optimal (in our class). In the second part of the article, we identify a second set of conditions, under which we prove the existence and strong consistency of a sequence of estimators of $\beta$, defined as roots of estimation equations which are martingale transforms (in particular, roots of the sequence of asymptotically optimal equations).

## 1 Introduction

### 1.1 Background

Longitudinal data sets are frequently used in biostatistics, economics, as well as in educational or environmental studies, when the individual measurements are recorded over time. Since in most applications, the individual measurements are influenced by a set of explanatory variables (e.g., age, family income, years of post-secondary education, etc.), the standard approach to longitudinal data analysis is based on a marginal regression model with unknown parameter $\beta$, which controls the effects of the explanatory variables. A study of longitudinal data is technically more demanding than a classical study of cross-sectional data, since it might have to include a separate model for the (unknown) correlation structure within the individual measurements. Such a study usually involves higher costs, but has the advantage of providing the researcher with the opportunity to control the unmeasured heterogenity among the response variables. In the most complex longitudinal scenarios, the observations are unbalanced (i.e. there is a different number of observations for each individual), the observational times depend on the individual, and the presence of some random/fixed effects specific to each individual cannot be ignored. We refer the reader to the monographs [5], [11] and [23] for a comprehensive account on this subject.

[*]Corresponding author. University of Ottawa. Department of Mathematics and Statistics. 585 King Edward Avenue, Ottawa, Ontario, K1N 6N5, Canada. E-mail address: rbalan@uottawa.ca

[†]This paper is based on a portion of the second author's doctoral thesis.

[‡]The first author and the third author are supported by research grants from the Natural Sciences and Engineering Research Council of Canada.

There are many approaches which have been proposed in the literature for treating the unknown correlation/covariance matrix within the individual responses. We discuss very briefly the salient features of some of these approaches.

The estimation of a covariance function in the context of correlated data (in particular longitudinal data) is an important statistical problem, for which various solutions have been proposed in the literature, using non-parametric methods [26], penalized likelihood methods [13], and methods borrowed from functional data analysis (see [28], [29]). More recently, the authors of [8] (see also [25]) have proposed a semi-parametric random-effect model for unbalanced time-dependent longitudinal data:

$$y_i(t) = \mathbf{x}_i(t)^T \beta + \delta_i(t)^T \alpha(t) + \varepsilon_i(t), \quad t = t_{ij}, j = 1, \ldots, m_i. \qquad (1)$$

In this model, the response $y_i(t)$ of the $i$-th individual at time $t$ depends linearly on a $p$-dimensional covariate vector $\mathbf{x}_i(t)$ (through the value of a regression parameter $\beta$), and a $d$-dimensional random-effect vector $\delta_i(t)$ (through the value of a smooth function $\alpha(t)$). The covariance structure within the individual responses is given by an unknown function $\sigma^2(t) := \text{Var}(\varepsilon(t)|\mathbf{x}(t), \delta(t))$ and a family $\rho(s, t; \theta) := \text{Corr}(\varepsilon(s), \varepsilon(t))$, which depends on a parameter $\theta$. The joint estimation of the correlation matrix and the regression parameter is achieved by iterating between the estimation of $(\sigma^2(t), \theta)$ and $(\alpha(t), \beta)$, using a combination of non-parametric and parametric techniques.

Other models related to (1) have been considered by various authors. Without aiming to exhaust the list of contributions in this active area of research, we mention briefly [30] (which underlies the connection with survival analysis), and [19] (motivated by an application in educational studies). In [30], the function $\alpha(t)$ is replaced by a normal random vector $\mathbf{a}$ which models the unobserved subject-specific effects, and the authors achieve a joint estimation of various parameters describing the "marker process" $y(t)$ and the survival time $T$. On the other hand, the authors of [19] considered the random-effect model (for balanced data) $y_i = \mathbf{X}_i\beta + \mathbf{A}\delta_i + \varepsilon_i$, where $\mathbf{X}_i = (\mathbf{x}_{i1}, \ldots, \mathbf{x}_{im})^T$, $\mathbf{A}$ is an $m \times d$ matrix of time-invariant individual parameters, and $m$ is the number of observations per individual.

The most popular approach for treating the unknown correlation structure (within the individual measurements) has its origins in [17], and has been embraced very quickly by the scientific community at large. It is now implemented in many statistical software. This approach is based on the marginal model:

$$y_{ij} = \mu(\mathbf{x}_{ij}^T\beta) + \varepsilon_{ij}, \quad j = 1, \ldots, m_i. \tag{2}$$

The main idea is to replace the true correlation matrix of the $i$-th individual, by a correlation matrix $\mathbf{R}_i(\alpha)$, which depends on a parameter $\alpha$, and achieve a joint estimation of $(\alpha, \beta)$ by iterating between the estimation of $\alpha$ and $\beta$. The usual recipe for the estimation of $\alpha$ is based on the method of moments. For the estimation of $\beta$, the authors of [17] suggest a quasi-likelihood method, inspired by the appealing similarity between a generalized linear model and the marginal model (2). This involves solving for $\beta$ in the *generalized estimating equation* (GEE):

$$\mathbf{g}_n(\beta) := \sum_{i=1}^n \mathbf{D}_i(\beta)^T \mathbf{V}_i^{-1}(\beta, \alpha)(\mathbf{y}_i - \mu_i(\beta)) = 0, \tag{3}$$

where $\mathbf{D}_i(\beta) = \partial \mu_i(\beta)/\partial \beta^T$, $\mu_i(\beta) = (\mu_{i1}(\beta), \ldots, \mu_{im_i}(\beta))^T$ and $\mu_{ij}(\beta) = E(y_{ij}) = \mu(\mathbf{x}_{ij}^T\beta)$. In this equation, $\mathbf{V}_i(\beta, \alpha) = \mathbf{A}_i(\beta)^{1/2}\mathbf{R}_i(\alpha)\mathbf{A}_i(\beta)^{1/2}$ is a "working covariance" matrix, which is obtained from the matrix $\mathbf{R}_i(\alpha)$ and the diagonal matrix $\mathbf{A}_i(\beta)$, built from the marginal variances $\sigma_{ij}^2(\beta) = \text{Var}(y_{ij}) = \mu'(\mathbf{x}_{ij}^T\beta)$. We refer the reader to [27] for an extensive theoretical study of various asymptotic properties of a different GEE estimator, as well as to [2] and [21] for similar studies of GEE's.

The appealing feature of the approach proposed in [17] is that the estimation of $\beta$ is "derived without specifying the joint distribution of the a subject's observations". Its drawback is that it requires a *correct* modelisation of the true correlation matrix from the very beginning.

One solution to this problem has been proposed recently in [14], using the theory of optimal parameter estimation, initiated by [9] and described at length in [12]. The model considered in [14] is written in semi-parametric form:

$$y_{ij} = g_j(\mathbf{X}_i, \beta) + \varepsilon_{ij}, \quad j = 1, , \ldots, m_i,$$

the term $\mu_{ij}(\beta) = g_j(\mathbf{X}_i, \beta)$ incorporating the subject-specific random effects. The authors of [14] propose an iterative procedure for the joint estimation of $\beta$ and the true covariance matrix $\Sigma_i$. The resulting *iterative estimating equations* (IEE) algorithm alternates between the estimation of $\beta$, and the method of moments estimation of $\Sigma_i$, converges exponentially fast, and yields consistent estimates for both $\beta$ and $\Sigma_i$. Most importantly, the IEE estimators for $\beta$ produced by this algorithm, are (asymptotically) as efficient as the optimal GEE estimators, defined as solutions of the equations:

$$\sum_{i=1}^{n} \mathbf{D}_i(\beta)^T \Sigma_i^{-1}(\mathbf{y}_i - \mu_i(\beta)) = 0, \quad n \geq 1. \tag{4}$$

## 1.2 Our contribution

In the present article, we consider a marginal model for balanced data, similar to (2), in which the covariates are non-random and there are no subject-specific random-effects. To simplify the presentation, we assume that the observational times are equally spaced and do not depend on the subjects. Our goals are to identify an asymptotically optimal equation (within an appropriate class), in the sense that it produces an estimator $\hat{\beta}_n$ which has minimum variance.

The model that we consider allows for some degree of dependence *among* the individual responses. In particular, our model includes the case in which the individual responses are independent. More precisely, we assume that the conditional mean and variance of the response $y_{ij}$ of individual $i$ at time $j$, given the "previous" responses $\mathbf{y}_1, \ldots, \mathbf{y}_{i-1}$ does *not* depend on these responses, and can be directly expressed using some explanatory variables $\mathbf{x}_{ij}$ (through the regression parameter $\beta$) and a link function $\mu$.

This theoretical relaxation has been inspired by the classical work [16] in the case of linear regression, in which the residuals are not necessarily independent, but form a martingale difference sequence. With these residuals as building blocks, we construct a class of estimating equations, which form a transform martingale family (see subsection 13.1.2 of [12]). The choice of this class leads us to an appropriate -and very general- class of estimating equations, in which a sequence of asymptotically optimal equations can be found. Although we follow the approach and use the techniques developed in [12], our application does not seem to have been covered in [12].

After the appropriate class of estimating equations has been selected (see p. 200 of [12]), the initial step of our investigation is to find the "optimal" equation (within this class). The method that we use for achieving this goal is different from the one of [14], since we allow that both the true conditional covariance matrix $\Sigma_i$ and the true correlation matrix $\overline{\mathbf{R}}_i$ depend on the parameter $\beta$. We find that the optimal estimating equation is:

$$\overline{\mathbf{g}}_n(\beta) := \sum_{i=1}^{n} \mathbf{D}_i(\beta)^T \Sigma_i(\beta)^{-1}(\mathbf{y}_i - \mu_i(\beta)) = 0, \tag{5}$$

which is significantly different from equation (4): even in the case of linear regression (i.e when $\mu(x) = x$), equation (5) cannot be solved explicitly, whereas (4) has a closed form solution, given by $\hat{\beta}_n = (\sum_{i=1}^{n} \mathbf{X}_i^T \Sigma_i^{-1} \mathbf{X}_i)^{-1} \sum_{i=1}^{n} \mathbf{X}_i^T \Sigma_i^{-1} \mathbf{y}_i$.

Next, we are interested in identifying an asymptotically optimal estimating equation, by means of a comparison with the optimal equation. To do this, we introduce a random matrix $\mathcal{R}_{i-1}^*(\beta)$ as an estimator of the true conditional correlation matrix $\overline{\mathbf{R}}_i(\beta)$. Under some conditions on the matrix $\mathcal{R}^*_{i-1}(\beta)$, we find that the following sequence of equations is asymptotically optimal:

$$\mathbf{g}^*_n(\beta) := \sum_{i=1}^n \mathbf{D}_i(\beta)^T \mathbf{V}^*_i(\beta)^{-1}(\mathbf{y}_i - \mu_i(\beta)) = 0, \quad n \geq 1 \qquad (6)$$

where $\mathbf{V}^*_i(\beta) = \mathbf{A}_i(\beta)^{1/2} \mathcal{R}^*_{i-1}(\beta) \mathbf{A}_i(\beta)^{1/2}$. A particular example of random matrix $\mathcal{R}^*_n(\beta)$, which can be used when the correlation structure is the same for all individuals, is:

$$\mathcal{R}^*_n(\beta) := \frac{1}{n} \sum_{i=1}^n \mathbf{A}_i(\beta)^{-1/2}(\mathbf{y}_i - \mu_i(\beta))(\mathbf{y}_i - \mu_i(\beta))^T \mathbf{A}_i(\beta)^{-1/2}. \qquad (7)$$

Computationally, solving this new equation may require more effort than (3) or (4), but it has the advantage of taking into account the correlation structure embedded in the data, without considering an additional model for this structure. An additional level of technical difficulty comes from the substitution of a random correlation matrix in (6), which renders the assumption of independence of residuals useless for the purpose of asymptotic analysis. These considerations have lead us to considering a class of transform martingales, which is required even in the case of independent individuals.

Most of the GEE literature deals with weak consistency of estimators. A result regarding the more difficult problem of strong convergence of the GEE estimators appears in [27]. The fortuitous choice of the above-mentioned class of estimating equations enables us to tackle this topic within a martingale framework and allows us to avoid using complicated approximation techniques (see Appendix A2, [2]). In the second part of the article (Section 4), we complete our analysis by specifying the sufficient conditions under which equation (6) can be solved, yielding a sequence of strongly consistent estimators of $\beta$. In addition to the previous example, this analysis applies also to the case when

$$\mathcal{R}^*_n(\beta) := \widetilde{\mathcal{R}}_n = \frac{1}{n} \sum_{i=1}^n \mathbf{A}_i(\tilde{\beta}_n)^{-1/2} \varepsilon_i(\tilde{\beta}_n) \varepsilon_i(\tilde{\beta}_n)^T \mathbf{A}_i(\tilde{\beta}_n)^{-1/2}. \qquad (8)$$

Here $\{\tilde{\beta}_n\}_n$ is a strongly consistent sequence of estimators, defined as roots of the "working independence" equation:

$$\mathbf{g}^{\text{indep}}_n(\beta) := \sum_{i=1}^n \mathbf{X}_i^T(\mathbf{y}_i - \mu_i(\beta)) = 0. \qquad (9)$$

This article is organized as follows. In Section 2, we introduce the general framework. In Section 3, we first prove that equation (5) is optimal in a certain class of estimating equations; then, we show that the sequence (6) of estimating equations is asymptotically optimal within the same class (Theorem 3.9). The main result of Section 4 (Theorem 4.13) identifies the conditions under which there exists a strongly consistent sequence $\{\hat{\beta}_n\}_n$ of estimators of $\beta$, defined as roots of equations (6).

Theorem 3.9 applies only to equations (6) corresponding to the sequence $\{\mathcal{R}_n^*(\beta)\}_n$ given by (7). On the other hand, Theorem 4.13 applies to equations (3) and (9) (although these equations are not asymptotically optimal), as well as equations (6) corresponding to sequences $\{\mathcal{R}_n^*(\beta)\}_n$ given by (7) or (8).

Some technical proofs which are needed in Section 4 are included in the appendices: Appendix A gives the formula for the calculation of the derivative of $\mathbf{g}_n^*(\beta)$, while Appendices B-E contain the proofs of some technical results. For these proofs, we use techniques which appear in [2] and [27], and results on strong convergence of martingales.

## 2 The Model Assumptions and Notation

We first specify the matrix notation, that we employ in the present article (see [22]). If $\lambda$ is a $p \times 1$ vector, we denote by $\|\lambda\|$ its Euclidean norm. If $\mathbf{A}$ is a $p \times p$ matrix, we denote with $\|\mathbf{A}\| = \sup_{\|\lambda\|=1} \|\mathbf{A}\lambda\|$ its operator norm, and with $\||A\|| = \sup_{\|\lambda\|=1} |\lambda^T \mathbf{A}\lambda|$ its spectral radius. If $\mathbf{A}$ is symmetric, then $\||\mathbf{A}\|| = \|\mathbf{A}\|$. We denote by $\det(\mathbf{A})$ the determinant of $\mathbf{A}$, and by $\text{tr}(\mathbf{A})$ the trace of $\mathbf{A}$. If $\mathbf{A}$ is a symmetric matrix, we denote by $\lambda_{\min}(\mathbf{A})$ and $\lambda_{\max}(\mathbf{A})$ its minimum eigenvalue, respectively its maximum eigenvalue. For any matrix $\mathbf{A}$, $\| \mathbf{A} \| = \{\lambda_{\max}(\mathbf{A}^T \mathbf{A})\}^{1/2}$. We let $\mathbf{A}^{1/2}$ be the symmetric square root of a positive definite matrix $\mathbf{A}$ and $\mathbf{A}^{-1/2} = (\mathbf{A}^{1/2})^{-1}$. Finally, we use the matrix notation $\mathbf{A} \leq \mathbf{B}$ if $\mathbf{B} - \mathbf{A}$ is non-negative definite, i.e. $\lambda^T \mathbf{A} \lambda \leq \lambda^T \mathbf{B} \lambda$ for any $p \times 1$ vector $\lambda$.

Finally, in this article, we denote by $C$ a generic constant which does not depend on $n$ and $\beta$, but is different from case to to case.

We now introduce the model assumptions and the estimating equation, which is the focus of investigation in the present article.

For each $i \geq 1$, let $\mathbf{y}_i = (y_{i1}, \ldots, y_{im})^T$, be the response variable of individual $i$, where $y_{ij}$ represents the response of individual $i$ at time $j$, and $m$ is a fixed time horizon, which is the same for all the individuals in the study. Clearly, the variables $(y_{ij})_{1 \leq j \leq m}$ display a non-trivial correlation structure, which, in the main application that we have in mind, is assumed to be the same for all individuals.

As in a classical regression problem, each outcome variable $y_{ij}$ is thought to have been influenced by a set of explanatory variables, whose values are given by a $p$-dimensional vector $\mathbf{x}_{ij}$. The following example illustrates the complexity of such a study.

One of our assumptions is that the explanatory variables $\mathbf{x}_{ij}$ are non-random and the response variables $(y_i)_{i \geq 1}$ are defined on a common probability space $(\Omega, \mathcal{F}, P_\beta)$. The uncertainty in this model is represented by the probability measure $P_\beta$, which depends on the unknown parameter $\beta \in \mathcal{T}$, where $\mathcal{T}$ is an open set in $\mathbb{R}^d$. This is a standard assumption in the theory of statistical inference. Another usual assumption encountered in the literature is that the fact that the variables $(y_i)_{i \geq 1}$ are independent under $P_\beta$, for any value of $\beta \in \mathcal{T}$.

Our model assumptions are: for each $\beta$, and for any $i \leq n, j \leq m$, we have:

$$\begin{aligned} E_\beta(y_{ij}|\mathcal{F}_{i-1}) &= \mu(\mathbf{x}_{ij}^T \beta) := \mu_{ij}(\beta) \\ \mathrm{Var}_\beta(y_{ij}|\mathcal{F}_{i-1}) &= \phi \mu'(\mathbf{x}_{ij}^T \beta) := \phi \sigma_{ij}^2(\beta), \end{aligned}$$

where $\mathcal{F}_{i-1}$ denotes the $\sigma$-field containing all the information about the variables $y_1, \ldots y_{i-1}$, $\mu$ is an arbitrary differentiable function with positive derivative, and $E_\beta(\cdot|\mathcal{F}_{i-1}), \mathrm{Var}_\beta(\cdot|\mathcal{F}_{i-1})$ denote the conditional expectation, respectively the conditional variance with respect to $P_\beta$. Here $\phi$ is a nuisance parameter; in what follows, we assume that $\phi = 1$.

Here are the most commonly used link functions $\mu$:
1. in the linear regression, $\mu(y) = y$;
2. in the log regression for count data, $\mu(y) = \exp(y)$;
3. in the logistic regression for binary data, $\mu(y) = \exp(y)/[1 + \exp(y)]$;
4. in the probit regression for binary data, $\mu(y) = \Phi(y)$, where $\Phi$ is the standard normal distribution function; we have $\dot{\Phi}(y) = (2\pi)^{-1/2} \exp(-y^2/2)$.

By definition, $(y_{ij} - \mu_{ij}(\beta))_{i \geq 1}$ is a martingale difference sequence, with respect to $P_\beta$, for any $j \leq m$.

We let $\mu_i(\beta) = (\mu_{i1}(\beta), \ldots, \mu_{im}(\beta))^T$, and $\mathbf{A}_i(\beta)$ be the diagonal matrix with entries $\sigma_{i1}^2(\beta), \ldots, \sigma_{im}^2(\beta)$.

Let $\Sigma_i^{(c)}(\beta)$ be the conditional covariance matrix of $\mathbf{y}_i$ given $\mathcal{F}_{i-1}$, with respect to $P_\beta$, whose elements are:

$$v_{i,jk}^{(c)}(\beta) := E_\beta[(y_{ij} - \mu_{ij}(\beta))(y_{ik} - \mu_{ik}(\beta))|\mathcal{F}_{i-1}], \quad 1 \leq j,k \leq m.$$

In matrix notation, we write

$$\Sigma_i^{(c)}(\beta) = E_\beta[(\mathbf{y}_i - \mu_i(\beta))(\mathbf{y}_i - \mu_i(\beta))^T|\mathcal{F}_{i-1}].$$

The matrix $\Sigma_i^{(c)}(\beta)$ has non-random elements $\sigma_{i1}^2(\beta), \ldots, \sigma_{im}^2(\beta)$ on the diagonal, but possibly random elements off the diagonal.

Some information about the dependence structure (with respect to $P_\beta$) *within* the components of $\mathbf{y}_i$ is contained in its conditional correlation matrix $\overline{\mathbf{R}}_i^{(c)}(\beta)$ given $\mathcal{F}_{i-1}$, whose elements are:

$$\overline{r}_{i,jk}^{(c)}(\beta) := \frac{v_{i,jk}^{(c)}(\beta)}{\sigma_{ij}(\beta)\sigma_{ik}(\beta)}, \quad 1 \leq j,k \leq m.$$

Note that $|\overline{r}_{i,jk}^{(c)}(\beta)| \leq 1$ $P_\beta$-a.s. and $\overline{r}_{i,jj}^{(c)}(\beta) = 1$ for any $\beta$. In matrix notation,

$$\Sigma_i^{(c)}(\beta) = \mathbf{A}_i(\beta)^{1/2} \overline{\mathbf{R}}_i^{(c)}(\beta) \mathbf{A}_i(\beta)^{1/2}. \tag{10}$$

Since $\partial \mu_{ij}(\beta)/\partial \beta^T = \sigma_{ij}^2(\beta) \mathbf{x}_{ij}^T$, which in matrix notation becomes:

$$\mathbf{D}_i(\beta) := \frac{\partial \mu_i(\beta)}{\partial \beta^T} = \mathbf{A}_i(\beta) \mathbf{X}_i,$$

where $\mathbf{X_i} = (\mathbf{x}_{i1}, \ldots, \mathbf{x}_{im})^T$ is an $m \times p$ matrix.

In the present article, we consider estimating equations of the form:

$$\mathbf{g}_n^*(\beta) := \sum_{i=1}^n \mathbf{D}_i(\beta)^T \mathbf{V}_i^*(\beta)^{-1}(\mathbf{y}_i - \mu_i(\beta)) = 0, \quad n \geq 1, \qquad (11)$$

where

$$\mathbf{V}_i^*(\beta) = \mathbf{A}_i(\beta)^{1/2} \mathcal{R}_{i-1}^*(\beta) \mathbf{A}_i(\beta)^{1/2} \qquad (12)$$

and $\{\mathcal{R}_n^*(\beta)\}_n$ is a sequence of random matrices, which satisfy the following conditions:

(A)  $\mathcal{R}_n^*(\beta)$ is positive-definite and continuously differentiable

(B)  the entries of $\mathcal{R}_n^*(\beta)$ are $\mathcal{F}_n$-measurable, for all $n \geq 1, \beta \in \mathcal{T}$.

One may think of the matrix $\mathcal{R}_{i-1}^*(\beta)$ as an approximation of the conditional correlation matrix $\overline{\mathbf{R}}_i^{(c)}(\beta)$, and hence of the matrix $\mathbf{V}_i^*(\beta)$ as an approximation of $\Sigma_i^{(c)}(\beta)$, due to (10) and (12). The fact that we consider $\mathcal{R}_{i-1}^*(\beta)$, instead of $\mathcal{R}_i^*(\beta)$, as an approximation for $\overline{\mathbf{R}}_i^{(c)}(\beta)$, guarantees that the function $\mathbf{g}_n^*(\beta)$ is a martingale.

The family $\{\mathbf{g}_n^*(\beta)\}_n$ is a transform martingale. (This family is a martingale with respect to $P_\beta$, if the entries of $\mathcal{R}_{i-1}^*(\beta)^{-1}(\mathbf{y}_i - \mu_i(\beta)$ are $P_\beta$-integrable.)

We now present several examples of estimating equations of the form (19).

**Example 2.1** The "working independence" estimating equations:

$$\mathbf{g}_n^{\text{indep}}(\beta) := \sum_{i=1}^n \mathbf{X}_i^T(\mathbf{y}_i - \mu_i(\beta)), \quad n \geq 1, \qquad (13)$$

constitute a particular case of (11), with $\mathcal{R}_{i-1}^*(\beta) = \mathbf{I}$ for all $i \geq 1$.

**Example 2.2** The "generalized estimating equations" (GEE) (3) studied in [27] can be written as:

$$\mathbf{g}_n^{\text{GEE}}(\beta) := \sum_{i=1}^n \mathbf{X}_i^T \mathbf{A}_i(\beta)^{1/2} \mathbf{R}_i(\alpha)^{-1} \mathbf{A}_i(\beta)^{-1/2}(\mathbf{y}_i - \mu_i(\beta)), \quad n \geq 1. \qquad (14)$$

These equations are particular instances of (11). In this case, $\mathcal{R}_{i-1}^*(\beta) = \mathbf{R}_i(\alpha)$ for all $i \geq 1$, where $\mathbf{R}_i(\alpha)$ are some non-random positive-definite matrices, depending on a parameter $\alpha$.

For the next two examples, we assume that the conditional correlation matrix $\overline{\mathbf{R}}_i^{(c)}(\beta)$ is the same for all individuals, i.e.

$$\overline{\mathbf{R}}_i^{(c)}(\beta) = \overline{\mathbf{R}}^{(c)}(\beta), \quad \forall i \geq 1, \forall \beta \in \mathcal{T}. \qquad (15)$$

**Example 2.3** As in [20], let $\{\widetilde{\beta}_n\}_n$ be a sequence of consistent estimators of $\beta_0$, defined as roots of $\mathbf{g}_n^{\text{indep}}(\beta) = 0$. (Under some conditions, one can prove that the sequence $\{\widetilde{\beta}_n\}_n$ exists; see e.g. [6] or Remark 4.16. Here $\beta_0$ is fixed and represents the true value of the parameter.)

If (15) holds, then under the conditions of Theorem 1 of [2] (and using an argument similar to the one used in the proof of this theorem), one can show that the sequence $\{\widetilde{\mathcal{R}}_n\}_n$ defined by

$$\widetilde{\mathcal{R}}_n := \frac{1}{n}\sum_{i=1}^n \mathbf{A}_i(\widetilde{\beta}_n)^{-1/2}(\mathbf{y}_i - \mu_i(\widetilde{\beta}_n))(\mathbf{y}_i - \mu_i(\widetilde{\beta}_n))^T \mathbf{A}_i(\widetilde{\beta}_n)^{-1/2}, \quad (16)$$

approximates the matrix $\overline{\mathbf{R}}^{(c)}(\beta_0)$, i.e. $\widetilde{\mathcal{R}}_n - \overline{\mathbf{R}}^{(c)}(\beta_0) \to 0$, element-wise, $P_{\beta_0}$-a.s. The following *pseudo-likelihood equations* (PLE's)

$$\widetilde{\mathbf{g}}_n(\beta) = \sum_{i=1}^n \mathbf{X}_i^T \mathbf{A}_i(\beta)^{1/2} \widetilde{\mathcal{R}}_{i-1}^{-1} \mathbf{A}_i(\beta)^{-1/2}(\mathbf{y}_i - \mu_i(\beta)) = 0, \quad n \geq 1, \quad (17)$$

constitute a particular case of (11) with $\mathcal{R}_{i-1}^*(\beta) = \widetilde{\mathcal{R}}_{i-1}$ for all $\beta \in \mathcal{T}$. (Note that (17) is different than equation (4) considered in [2], which contains $\widetilde{\mathcal{R}}_n^{-1}$ in the middle, instead of $\widetilde{\mathcal{R}}_{i-1}^{-1}$.)

**Example 2.4** Suppose that (15) holds, and there exist some constants $C_\beta > 0$ and $\delta_\beta > 0$ such that

$$E_\beta \|\mathbf{A}_i(\beta)^{-1/2}(\mathbf{y}_i - \mu_i(\beta))\|^{2+\delta_\beta} \leq C_\beta, \quad \forall i \geq 1.$$

Using Lemma A.1 of [2], one can show that the sequence $\{\mathcal{R}_n^*(\beta)\}_n$, defined by

$$\mathcal{R}_n^*(\beta) := \frac{1}{n}\sum_{i=1}^n \mathbf{A}_i(\beta)^{-1/2}(\mathbf{y}_i - \mu_i(\beta))(\mathbf{y}_i - \mu_i(\beta))^T \mathbf{A}_i(\beta)^{-1/2}, \quad \beta \in \mathcal{T}, n \geq 1, \quad (18)$$

approximates the matrix $\overline{\mathbf{R}}^{(c)}(\beta)$, i.e. $\mathcal{R}_{n-1}^*(\beta) - \overline{\mathbf{R}}^{(c)}(\beta) \to 0$ element-wise, $P_\beta$-a.s and in $L^1(P_\beta)$, for all $\beta \in \mathcal{T}$. The sequence $\{\mathcal{R}_n^*(\beta)\}_n$ satisfies conditions $(A)$ and $(B)$. Equation (11) can be written as:

$$\mathbf{g}_n^*(\beta) = \sum_{i=1}^n \mathbf{X}_i^T \mathbf{A}_i(\beta)^{1/2} \mathcal{R}_{i-1}^*(\beta)^{-1} \mathbf{A}_i(\beta)^{-1/2}(\mathbf{y}_i - \mu_i(\beta)) = 0, \quad n \geq 1. \quad (19)$$

(Similar estimating equations, which were not transform martingales were studied in [21].)

We consider the following sequence of estimating equations:

$$\overline{\mathbf{g}}_n(\beta) := \sum_{i=1}^n \mathbf{D}_i(\beta)^T \Sigma_i^{(c)}(\beta)(\mathbf{y}_i - \mu_i(\beta)) = 0, \quad n \geq 1,$$

which can be written as:

$$\overline{\mathbf{g}}_n(\beta) = \sum_{i=1}^{n} \mathbf{X}_i^T \mathbf{A}_i(\beta)^{1/2} \overline{\mathbf{R}}_i^{(c)}(\beta)^{-1} \mathbf{A}_i(\beta)^{-1/2}(\mathbf{y}_i - \mu_i(\beta)) = 0, \quad n \geq 1.$$

We have:

$$\mathbf{M}_n^*(\beta) := \mathrm{Cov}_\beta[\mathbf{g}_n^*(\beta)] = \sum_{i=1}^{n} \mathbf{X}_i^T \mathbf{A}_i(\beta)^{1/2} \overline{\mathbf{E}}_{i-1}^*(\beta) \mathbf{A}_i(\beta)^{1/2} \mathbf{X}_i \quad (20)$$

$$\overline{\mathbf{M}}_n(\beta) := \mathrm{Cov}_\beta[\overline{\mathbf{g}}_n(\beta)] = \sum_{i=1}^{n} \mathbf{X}_i^T \mathbf{A}_i(\beta)^{1/2} \overline{\mathbf{E}}_i(\beta) \mathbf{A}_i(\beta)^{1/2} \mathbf{X}_i. \quad (21)$$

Here $\overline{\mathbf{E}}_{i-1}^*(\beta) := E_\beta[\mathcal{R}_{i-1}^*(\beta)^{-1} \overline{\mathbf{R}}_i^{(c)}(\beta) \mathcal{R}_{i-1}^*(\beta)^{-1}]$ and $\overline{\mathbf{E}}_i(\beta) := E_\beta[\overline{\mathbf{R}}_i^{(c)}(\beta)^{-1}]$.

## 3 Optimal Estimating Equation

Following the approach of [12], we introduce a general class $\mathcal{H}_n$ of estimating functions (which accommodate our model), and the concept of optimal estimating equation in this class. As a preliminary step, we show that the estimating function $\overline{\mathbf{g}}_n(\beta)$ is optimal within this class. The main result of this section identifies a set of conditions for the approximation matrices $\{\mathcal{R}_n^*(\beta)\}_{n\geq 1}$, under which the sequence $\{\mathbf{g}_n^*(\beta)\}_{n\geq 1}$ of estimating equations is "asymptotically optimal" within $\{\mathcal{H}_n\}_{n\geq 1}$.

For each $n \geq 1$, we consider the following class of estimating functions:

$$\mathcal{H}_n = \{\mathbf{q}_n(\beta) = \sum_{i=1}^{n} \mathbf{C}_i(\beta)(\mathbf{y}_i - \mu_i(\beta)), \beta \in \mathcal{T}\},$$

where $\mathbf{C}_i(\beta)$ is a $p \times m$ random matrix, whose elements are $\mathcal{F}_{i-1}$-measurable and continuously differentiable (with respect to $\beta$), for all $i \geq 1$. Moreover, if $c_{i,uj}(\beta)$ denotes the $(u,j)$-element of $\mathbf{C}_i(\beta)$, we assume that: for any $\beta \in \mathcal{T}$, $i \geq 1$, $1 \leq u, v \leq p$ and $1 \leq j, k \leq m$

$$E_\beta |c_{i,uj}(\beta)| < \infty, \quad E_\beta\left[\frac{\partial c_{i,uj}(\beta)}{\partial \beta_v}(\mathbf{y}_{ij} - \mu_{ij}(\beta))\right] < \infty,$$

$$E_\beta[c_{i,uj}(\beta) v_{i,jk}^{(c)}(\beta) c_{i,vk}(\beta)] < \infty.$$

For each function $\mathbf{q}_n(\beta) \in \mathcal{H}_n$, we introduce the following matrix:

$$\mathcal{E}[\mathbf{q}_n(\beta)] := \left\{ E_\beta\left[\frac{\partial \mathbf{q}_n(\beta)}{\partial \beta^T}\right]\right\}^T \{\mathrm{Cov}_\beta[\mathbf{q}_n(\beta)]\}^{-1} E_\beta\left[\frac{\partial \mathbf{q}_n(\beta)}{\partial \beta^T}\right].$$

**Remark 3.1** Note that $\overline{\mathbf{g}}_n(\beta)$ is an element of the class $\mathcal{H}_n$. Another element of $\mathcal{H}_n$ is the GEE function $\mathbf{g}_n(\beta)$ of [27], given by (14).

**Remark 3.2** The function $\mathbf{g}_n^{\text{indep}}(\beta)$, given by (13) is also an element of $\mathcal{H}_n$. For this function, we have:

$$\mathbf{H}_n^{\text{indep}}(\beta) := -E_\beta \left[ \frac{\partial \mathbf{g}_n^{\text{indep}}(\beta)}{\partial \beta^T} \right] = \sum_{i=1}^n \mathbf{X}_i^T \mathbf{A}_i(\beta) \mathbf{X}_i$$

$$\mathbf{M}_n^{\text{indep}}(\beta) := \text{Cov}_\beta[\mathbf{g}_n^{\text{indep}}(\beta)] = \sum_{i=1}^n \mathbf{X}_i^T \mathbf{A}_i(\beta)^{1/2} \overline{\mathbf{R}}_i(\beta) \mathbf{A}_i(\beta)^{1/2} \mathbf{X}_i,$$

where $\overline{\mathbf{R}}_i(\beta)$ is the true (unconditional) correlation matrix of the $i$-th individual.

The function $\mathbf{g}_n^{\text{indep}}(\beta)$ can be viewed as a score function, in a model in which $\overline{\mathbf{R}}_i(\beta) = \mathbf{I}$. (Recall that a score function is the derivative of a log-likelihood function.) To see this, suppose that there exists a function $a$ such that $a' = \mu$. Then $\mathbf{g}_n^{\text{indep}}(\beta) = \partial l_n(\beta)/\partial \beta^T$, where $l_n(\beta) = \sum_{i=1}^n [\beta^T \mathbf{X}_i^T y_i - \sum_{j=1}^m a(\mathbf{x}_{ij}^T \beta)]$.

**Remark 3.3** Let $\mathbf{g}_n(\beta)$ be the GEE function, given by (14). Let $\mathbf{M}_n := \text{Cov}_{\beta_0}[\mathbf{g}_n(\beta_0)]$ and $\mathbf{H}_n := E_{\beta_0}[\mathcal{D}_n(\beta_0))]$, where $\mathcal{D}_n(\beta) = -\partial \mathbf{g}_n(\beta)/\partial \beta^T$. (Here $\beta_0 \in \mathcal{T}$ is fixed and represents the "true" value of the parameter.) If $\{\hat{\beta}_n\}_n$ is a sequence of weakly consistent estimators of $\beta_0$, then Theorem 4 of [27] says that

$$\mathbf{M}_n^{-1/2} \mathbf{H}_n(\hat{\beta}_n - \beta_0) \to N(0, \mathbf{I}),$$

in distribution (under $P_{\beta_0}$). Therefore, in order to obtain an asymptotic confidence interval for $\beta_0$ of minimal length, one needs to *maximize* (in the sense of the non-negative definiteness order) the matrix $\mathbf{H}_n \mathbf{M}_n^{-1} \mathbf{H}_n = \mathcal{E}[\mathbf{g}_n(\beta_0)]$. This gives a first motivation for Definition 3.5.

**Remark 3.4** The matrix $\mathcal{E}[\mathbf{q}_n(\beta)]$ can be viewed as a generalization of Fisher information matrix. To see this, recall that if $\mathbf{s}_n(\beta)$ is a score function in the class $\mathcal{H}_n$, then $\mathcal{E}[\mathbf{s}_n(\beta)]$ coincides with Fisher information matrix:

$$\mathcal{E}[\mathbf{s}_n(\beta)] = \text{Cov}_\beta[\mathbf{s}_n(\beta)] = -E_\beta \left[ \frac{\partial \mathbf{s}_n(\beta)}{\partial \beta^T} \right].$$

By the Cramér-Rao inequality (see e.g. Theorem 7.3.10, [3]), the best unbiased estimator $W = W(\mathbf{y}_1, \ldots, \mathbf{y}_n)$ of $\beta$ is the one which attains the Cramér-Rao lower bound, i.e. for which $\text{Cov}_\beta[W] = \{\mathcal{E}[\mathbf{s}_n(\beta)]\}^{-1}$. Among those estimators, the one with minimum variance is the one for which Fisher information matrix is maximal. This provides another motivation for Definition 3.5.

We are now ready to introduce the concept of optimal estimating function in the class $\mathcal{H}_n$ (see Definition 2.1, [12]).

**Definition 3.5** *We say that an estimating function $\mathbf{q}_n^*(\beta) \in \mathcal{H}_n$ is **optimal** (or **quasi-score**) within the class $\mathcal{H}_n$, if for any $\mathbf{q}_n(\beta) \in \mathcal{H}_n$ and for any $\beta \in \mathcal{T}$*

$$\mathcal{E}[\mathbf{q}_n^*(\beta)] - \mathcal{E}[\mathbf{q}_n(\beta)] \quad \text{is nonnegative-definite.}$$

The following result lies at the origin of our developments.

**Proposition 3.6** *The function $\overline{\mathbf{g}}_n(\beta)$ is a quasi-score within the class $\mathcal{H}_n$.*

**Proof:** Using Theorem 2.1, [12], it suffices to show that for any $\mathbf{q}_n \in \mathcal{H}_n$,

$$E_\beta[\mathbf{q}_n(\beta)\overline{\mathbf{g}}_n(\beta)^T] = -E_\beta\left[\frac{\partial \mathbf{q}_n(\beta)}{\partial \beta^T}\right], \quad \forall \beta \in \mathcal{T}. \tag{22}$$

First, we treat the left-hand side of (22). Note that $\overline{\mathbf{g}}_n(\beta) = \sum_{i=1}^n \overline{\mathbf{C}}_i(\beta)\varepsilon_i(\beta)$, where

$$\overline{\mathbf{C}}_i(\beta) := \mathbf{X}_i^T \mathbf{A}_i(\beta)^{1/2} \overline{\mathbf{R}}_i^{(c)}(\beta)^{-1} \mathbf{A}_i(\beta)^{-1/2} = \mathbf{X}_i^T \mathbf{A}_i(\beta) \Sigma_i^{(c)}(\beta)^{-1}.$$

Using the fact that $E_\beta(\mathbf{y}_i|\mathcal{F}_{i-1}) = \mu_i(\beta)$, we obtain:

$$\begin{aligned}
E_\beta[\mathbf{q}_n(\beta)\overline{\mathbf{g}}_n(\beta)^T] &= \sum_{i=1}^n E_\beta\{\mathbf{C}_i(\beta) E_\beta[(\mathbf{y}_i - \mu_i(\beta))(\mathbf{y}_i - \mu_i(\beta))^T | \mathcal{F}_{i-1}] \overline{\mathbf{C}}_i(\beta)^T\} \\
&= \sum_{i=1}^n E_\beta[\mathbf{C}_i(\beta) \Sigma_i^{(c)}(\beta) \overline{\mathbf{C}}_i(\beta)^T] \\
&= \sum_{i=1}^n E_\beta[\mathbf{C}_i(\beta) \mathbf{A}_i(\beta) \mathbf{X}_i]. \tag{23}
\end{aligned}$$

Next, we treat the right-hand side of (22). If we denote by $\mathbf{c}_{i1}(\beta), \ldots, \mathbf{c}_{im}(\beta)$ the columns of $\mathbf{C}_i(\beta)$, then $\mathbf{q}_n(\beta) = \sum_{i=1}^n \sum_{j=1}^m \mathbf{c}_{ij}(\beta)(y_{ij} - \mu_{ij}(\beta))$. Using the chain rule, the fact that $\mathbf{c}_{ij}(\beta)$ is $\mathcal{F}_{i-1}$-measurable, and $E_\beta(y_{ij}|\mathcal{F}_{i-1}) = \mu_{ij}(\beta)$, we have:

$$\begin{aligned}
E_\beta\left[\frac{\partial \mathbf{q}_n(\beta)}{\partial \beta^T}\right] &= \sum_{i=1}^n \sum_{j=1}^m \left\{ E_\beta\left[\frac{\partial \mathbf{c}_{ij}(\beta)}{\partial \beta^T}(y_{ij} - \mu_{ij}(\beta))\right] - E_\beta\left[\mathbf{c}_{ij}(\beta)\frac{\partial \mu_{ij}(\beta)}{\partial \beta^T}\right]\right\} \\
&= \sum_{i=1}^n \sum_{j=1}^m \left\{ E_\beta\left[\frac{\partial \mathbf{c}_{ij}(\beta)}{\partial \beta^T} E_\beta(y_{ij} - \mu_{ij}(\beta))|\mathcal{F}_{i-1})\right] - \right. \\
&\qquad E_\beta\left[\mathbf{c}_{ij}(\beta)\sigma_{ij}^2(\beta)\mathbf{x}_{ij}^T\right]\} \\
&= -\sum_{i=1}^n \sum_{j=1}^m E_\beta\left[\mathbf{c}_{ij}(\beta)\sigma_{ij}^2(\beta)\mathbf{x}_{ij}^T\right] \\
&= -\sum_{i=1}^n E_\beta[\mathbf{C}_i(\beta)\mathbf{A}_i(\beta)\mathbf{X}_i]. \tag{24}
\end{aligned}$$

Relation (22) follows from (23) and (24) □

**Remark 3.7** By taking $\mathbf{q}_n(\beta) = \overline{\mathbf{g}}_n(\beta)$ in (22), we see that $\overline{\mathbf{g}}_n(\beta)$ has the property of a score function:

$$\overline{\mathbf{H}}_n(\beta) := -E_\beta\left[\frac{\partial \overline{\mathbf{g}}_n(\beta)}{\partial \beta^T}\right] = \text{Cov}_\beta[\overline{\mathbf{g}}_n(\beta)] = \overline{\mathbf{M}}_n(\beta). \tag{25}$$

Also, by taking $\mathbf{q}_n(\beta) = \mathbf{g}_n^*(\beta)$ in (22), we obtain that:

$$\mathbf{H}_n^*(\beta) := -E_\beta\left[\frac{\partial \mathbf{g}_n^*(\beta)}{\partial \beta^T}\right] = \sum_{i=1}^n \mathbf{X}_i^T \mathbf{A}_i(\beta)^{1/2} \mathbf{E}_{i-1}^*(\beta) \mathbf{A}_i(\beta)^{1/2} \mathbf{X}_i, \qquad (26)$$

where $\mathbf{E}_{i-1}^*(\beta) := E_\beta[\mathcal{R}_{i-1}^*(\beta)^{-1}]$. From here, we conclude that:

$$\mathcal{E}[\overline{\mathbf{g}}_n(\beta)] = \overline{\mathbf{M}}_n(\beta) \quad \text{and} \quad \mathcal{E}[\mathbf{g}_n^*(\beta)] = \mathbf{H}_n^*(\beta)\mathbf{M}_n^*(\beta)^{-1}\mathbf{H}_n^*(\beta). \qquad (27)$$

We note that the optimal function $\overline{\mathbf{g}}_n(\beta)$ depends on the unknown conditional correlation matrix $\overline{\mathbf{R}}_i^{(c)}(\beta)$, and therefore, cannot be used in practice. In the remaining part of this section, we circumvent this difficulty by replacing it with a consistent estimator $\mathcal{R}_{i-1}^*(\beta)$, proving that this procedure preserves the optimality of the equation, in the asymptotic sense.

As in [12], we consider now the "normalized" estimating function:

$$\mathbf{q}_n^{(\text{norm})}(\beta) := \left\{E_\beta\left[\frac{\partial \mathbf{q}_n(\beta)}{\partial \beta^T}\right]\right\}^{-1} \mathbf{q}_n(\beta),$$

for any $\mathbf{q}_n(\beta) \in \mathcal{H}_n$. Note that the covariance matrix of $\mathbf{q}_n^{(\text{norm})}(\beta)$ is:

$$\mathcal{I}[\mathbf{q}_n(\beta)] := E_\beta[\mathbf{q}_n^{(\text{norm})}(\beta)\mathbf{q}_n^{(\text{norm})}(\beta)^T] = \{\mathcal{E}[\mathbf{q}_n(\beta)]\}^{-1}. \qquad (28)$$

A sequence $\{\mathbf{q}_n^*(\beta)\}_{n\geq 1}$ of estimating functions is asymptotically optimal, within the collection $\{\mathcal{H}_n\}_{n\geq 1}$, if the corresponding matrix $\mathcal{I}[\mathbf{q}_n^*(\beta)]$ is *minimal* (in the sense of the non-negative definiteness order), when $n$ is large enough.

More precisely, we have the following definition: (see Definition 5.1, [12])

**Definition 3.8** *Let $\{\mathbf{q}_n^*(\beta)\}_{n\geq 1}$ be a sequence of estimating functions such that $\mathbf{q}_n^*(\beta) \in \mathcal{H}_n$ for all $n \geq 1$. We say that $\{\mathbf{q}_n^*(\beta)\}_{n\geq 1}$ is **asymptotically optimal** (or **asymptotic quasi-score**) within the collection $\{\mathcal{H}_n\}_{n\geq 1}$, if for any sequence $\{\mathbf{q}_n(\beta)\}_{n\geq 1}$ with $\mathbf{q}_n(\beta) \in \mathcal{H}_n$ for all $n \geq 1$, and for any $\beta \in \mathcal{T}$,*

$$\{\mathcal{I}_n^*(\beta)^{-1/2}\mathcal{I}_n(\beta)\mathcal{I}_n^*(\beta)^{-1/2} - \mathbf{I}\}_{n\geq 1} \text{ is asymptotically non-negative definite,}$$

*in the sense that, and for any $p \times 1$ vector $\lambda$ with $\|\lambda\| = 1$,*

$$\liminf_{n\to\infty} \lambda^T[\mathcal{I}_n^*(\beta)^{-1/2}\mathcal{I}_n(\beta)\mathcal{I}_n^*(\beta)^{-1/2} - \mathbf{I}]\lambda \geq 0.$$

*Here $\mathcal{I}_n^*(\beta) = \mathcal{I}[\mathbf{q}_n^*(\beta)]$ and $\mathcal{I}_n(\beta) = \mathcal{I}[\mathbf{q}_n(\beta)]$.*

(See Remark 5.3 of [12] for a motivation of the previous definition, and the proof of Proposition 5.4 of [12] for the rigorous meaning of the concept of "asymptotic non-negative definiteness" introduced above.)

Similarly to [2], we introduce the following assumption:

(H)    there exists a constant $C_\beta > 0$ such that $\lambda_{\min}[\overline{\mathbf{R}}_n^{(c)}(\beta)] \geq C_\beta, \forall n \geq 1$, $P_\beta$-almost surely, for all $\beta \in \mathcal{T}$.

(see condition $(H')$ on p. 528 of [2])

The following theorem is the main result of this section.

**Theorem 3.9** *Suppose that assumption $(H)$ holds and*

$$\lambda_{\min}[\mathbf{H}_n^{\text{indep}}(\beta)] \to \infty, \text{ for all } \beta \in \mathcal{T}. \tag{29}$$

*Let $\{\mathcal{R}_n^*(\beta)\}_{n\geq 1}$ be a sequence of random matrices which satisfy conditions (A), (B), as well as the following conditions:*

(C)  $\mathcal{R}_{n-1}^*(\beta) - \overline{\mathbf{R}}_n^{(c)}(\beta) \to 0$ *(element-wise), in probability $P_\beta$, for all $\beta \in \mathcal{T}$*

(R)  *there exists a constant $K_\beta > 0$ such that $\lambda_{\min}[\mathcal{R}_n^*(\beta)] \geq K_\beta$ for all $n \geq 1$, $P_\beta$-almost surely, for all $\beta \in \mathcal{T}$.*

*Then, the sequence $\{\mathbf{g}_n^*(\beta)\}_{n\geq 1}$ is an asymptotic quasi-score within the collection $\{\mathcal{H}_n\}_{n\geq 1}$.*

**Proof:** By Proposition 3.6 and Remark 5.2, [12], the sequence $\{\overline{\mathbf{g}}_n(\beta)\}_{n\geq 1}$ is an asymptotic quasi-score within the collection $\{\mathcal{H}_n\}_{n\geq 1}$. By invoking Proposition 5.5, [12], it suffices to show that the sequences $\{\mathbf{g}_n^*(\beta)\}_{n\geq 1}$ and $\{\overline{\mathbf{g}}_n(\beta)\}_{n\geq 1}$ are "asymptotically equivalent", in the sense that

$$\frac{\det \mathcal{I}[\overline{\mathbf{g}}_n(\beta)]}{\det \mathcal{I}[\mathbf{g}_n^*(\beta)]} \to 1, \quad \forall \beta \in \mathcal{T}.$$

By (28), this is equivalent to:

$$\frac{\det \mathcal{E}[\mathbf{g}_n^*(\beta)]}{\det \mathcal{E}[\overline{\mathbf{g}}_n(\beta)]} \to 1, \quad \forall \beta \in \mathcal{T},$$

which in turn, by (27), is equivalent to:

$$\frac{\det \mathbf{H}_n^*(\beta)^2}{\det [\overline{\mathbf{M}}_n(\beta)\mathbf{M}_n^*(\beta)]} \to 1, \quad \forall \beta \in \mathcal{T}.$$

Therefore, the proof will be complete, once we show that for any $\beta \in \mathcal{T}$

$$\frac{\det \mathbf{H}_n^*(\beta)}{\det \overline{\mathbf{M}}_n(\beta)} \to 1 \quad \text{and} \quad \frac{\det \mathbf{M}_n^*(\beta)}{\det \overline{\mathbf{M}}_n(\beta)} \to 1. \tag{30}$$

Recalling the definitions (20), (21) and (26) of $\mathbf{M}_n^*(\beta)$, $\overline{\mathbf{M}}_n(\beta)$ and $\mathbf{H}_n^*(\beta)$, we see that to prove (30), it suffices to compare $\overline{\mathbf{E}}_i(\beta) := E_\beta[\overline{\mathbf{R}}_i^{(c)}(\beta)^{-1}]$ with $\mathbf{E}_{i-1}^*(\beta) = E_\beta[\mathcal{R}_{i-1}^*(\beta)^{-1}]$ and $\overline{\mathbf{E}}_{i-1}^*(\beta) = E_\beta[\mathcal{R}_{i-1}^*(\beta)^{-1}\overline{\mathbf{R}}_i^{(c)}(\beta)\mathcal{R}_{i-1}^*(\beta)^{-1}]$.

Using $(H)$, $(C)$ and $(R)$, we claim that: (see below)

$$\begin{aligned}
\overline{\mathbf{E}}_i(\beta)^{-1}\mathbf{E}_{i-1}^*(\beta) &\to \mathbf{I} \quad \text{(elementwise)} & (31) \\
\overline{\mathbf{E}}_i(\beta)^{-1}\overline{\mathbf{E}}_{i-1}^*(\beta) &\to \mathbf{I} \quad \text{(elementwise)}. & (32)
\end{aligned}$$

We now proceed with the proof of the first convergence in (30), using (31); the second convergence follows by a similar argument, using (32). Let $\varepsilon \in (0,1)$ be arbitrary. By (31), there exists an integer $n_0$ (depending on $\varepsilon$ and $\beta$), such that

$$1 - \varepsilon \leq \lambda_{\min}[\overline{\mathbf{E}}_i(\beta)^{-1}\mathbf{E}_{i-1}^*(\beta)] \leq \lambda_{\max}[\overline{\mathbf{E}}_i(\beta)^{-1}\mathbf{E}_{i-1}^*(\beta)] \leq 1 + \varepsilon, \quad \forall i \geq n_0.$$

Therefore,

$$(1-\varepsilon)\overline{\mathbf{M}}_{n_0,n}(\beta) \leq \mathbf{H}_{n_0,n}^*(\beta) \leq (1+\varepsilon)\overline{\mathbf{M}}_{n_0,n}(\beta), \quad \forall n \geq n_0, \qquad (33)$$

where $\mathbf{H}_{n_0,n}^*(\beta) := \sum_{i=n_0}^{n} \mathbf{X}_i^T \mathbf{A}_i(\beta)^{1/2} \mathbf{E}_{i-1}^*(\beta) \mathbf{A}_i(\beta)^{1/2} \mathbf{X}_i$ and $\overline{\mathbf{M}}_{n_0,n}(\beta) := \sum_{i=n_0}^{n} \mathbf{X}_i^T \mathbf{A}_i(\beta)^{1/2} \overline{\mathbf{E}}_i(\beta) \mathbf{A}_i(\beta)^{1/2} \mathbf{X}_i$. Using the fact that the determinant is a non-decreasing function (with respect to the non-negative definiteness order), we obtain:

$$(1-\varepsilon)^p \leq \frac{\det \mathbf{H}_{n_0,n}^*(\beta)}{\det \overline{\mathbf{M}}_{n_0,n}(\beta)} \leq (1+\varepsilon)^p, \quad \forall n \geq n_0. \qquad (34)$$

Since $\overline{\mathbf{E}}_{i-1}(\beta) \geq m^{-1}\mathbf{I}$, it follows that $\overline{\mathbf{M}}_n(\beta) \geq m^{-1}\mathbf{H}_n^{\text{indep}}(\beta)$. By (29), $\lambda_{\min}[\overline{\mathbf{M}}_n(\beta)] \to \infty$, and therefore $\lambda_{\min}[\overline{\mathbf{M}}_{n_0,n}(\beta)] \to \infty$ as $n \to \infty$. Hence, there exists an integer $n_1 > n_0$ (depending on $\varepsilon$ and $\beta$) such that $\lambda_{\min}[\overline{\mathbf{M}}_{n_0,n}(\beta)] \geq \varepsilon^{-1}\lambda_{\max}[\overline{\mathbf{M}}_{n_0-1}(\beta)]$ for all $n \geq n_1$. Therefore, $\overline{\mathbf{M}}_{n_0-1}(\beta) \leq \varepsilon \overline{\mathbf{M}}_{n_0,n}(\beta)$ for all $n \geq n_1$, and

$$\overline{\mathbf{M}}_{n_0,n}(\beta) \leq \overline{\mathbf{M}}_n(\beta) \leq (1+\varepsilon)\overline{\mathbf{M}}_{n_0,n}(\beta), \quad \forall n \geq n_1.$$

From here, we conclude that:

$$\det \overline{\mathbf{M}}_{n_0,n}(\beta) \leq \det \overline{\mathbf{M}}_n(\beta) \leq (1+\varepsilon)^p \det \overline{\mathbf{M}}_{n_0,n}(\beta), \quad \forall n \geq n_1. \qquad (35)$$

This argument can be repeated for $\mathbf{H}_n^*(\beta)$, since $\lambda_{\min}[\mathbf{H}_{n_0,n}^*(\beta)] \to \infty$ as $n \to \infty$ (this is a consequence of (33)). We conclude that there exists an integer $n_2 > n_1$ (depending on $\varepsilon$ and $\beta$) such that

$$\det \mathbf{H}_{n_0,n}^*(\beta) \leq \det \mathbf{H}_n^*(\beta) \leq (1+\varepsilon)^p \det \mathbf{H}_{n_0,n}^*(\beta), \quad \forall n \geq n_2. \qquad (36)$$

From (35) and (36), we obtain:

$$\frac{1}{(1+\varepsilon)^p} \frac{\det \mathbf{H}_{n_0,n}^*(\beta)}{\det \overline{\mathbf{M}}_{n_0,n}(\beta)} \leq \frac{\det \mathbf{H}_n^*(\beta)}{\det \overline{\mathbf{M}}_n(\beta)} \leq (1+\varepsilon)^p \frac{\det \mathbf{H}_{n_0,n}^*(\beta)}{\det \overline{\mathbf{M}}_{n_0,n}(\beta)}, \quad \forall n \geq n_2. \qquad (37)$$

Finally, using (34) and (37), we obtain:

$$\left(\frac{1-\varepsilon}{1+\varepsilon}\right)^p \leq \frac{\det \overline{\mathbf{M}}_n(\beta)}{\det \mathbf{H}_n^*(\beta)} \leq (1+\varepsilon)^{2p}, \quad \forall n \geq n_2.$$

This concludes the proof of (30).

We now turn to the proof of (31). Suppose, by contradiction, that there exist some $\varepsilon_0 > 0$ and a subsequence $(i_n)_n$ for which

$$\|\overline{\mathbf{E}}_{i_n}(\beta)^{-1}\mathbf{E}^*_{i_n-1}(\beta) - \mathbf{I}\| > \varepsilon_0, \quad \forall n \geq 1. \tag{38}$$

Under condition (C), there exists a subsequence of $(i_n)_n$, which we denote by $(l_n)_n$, for which $\mathcal{R}^*_{l_n-1}(\beta) - \overline{\mathbf{R}}^{(c)}_{l_n}(\beta) \to 0$ (element-wise), $P_\beta$-a.s. Therefore, $\mathcal{R}^*_{l_n-1}(\beta)^{-1} - \overline{\mathbf{R}}^{(c)}_{l_n}(\beta)^{-1} \to 0$ (element-wise) $P_\beta$-a.s. By conditions $(R)$ and $(H)$, each of the elements of the matrices $\mathcal{R}^*_{l_n-1}(\beta)^{-1}$ and $\overline{\mathbf{R}}^{(c)}_{l_n}(\beta)^{-1}$ are bounded above by the constants $K_\beta^{-1}$, respectively $C_\beta^{-1}$, $P_\beta$-a.s., and hence this last convergence holds in $L^1(P_\beta)$ as well, i.e.

$$\mathbf{E}^*_{l_n-1}(\beta) - \overline{\mathbf{E}}_{l_n}(\beta) \to 0 \quad \text{(element-wise)}. \tag{39}$$

(Clearly, the element-wise convergence is equivalent to the convergence in norm.)

Note that $\overline{\mathbf{R}}^{(c)}_{l_n}(\beta)$ is an $m \times m$ positive-definite matrix, whose elements are bounded above by 1, $P_\beta$-a.s. Hence, $\overline{\mathbf{R}}^{(c)}_{l_n}(\beta) \leq m\mathbf{I}$, $P_\beta$-a.s. and $\overline{\mathbf{R}}^{(c)}_{l_n}(\beta)^{-1} \geq m^{-1}\mathbf{I}$, $P_\beta$-a.s. From here, we conclude that

$$\overline{\mathbf{E}}_{l_n}(\beta) = E_\beta[\overline{\mathbf{R}}^{(c)}_{l_n}(\beta)^{-1}] \geq m^{-1}\mathbf{I}, \quad \text{i.e.} \quad \lambda_{\min}[\overline{\mathbf{E}}_{l_n}(\beta)] \geq m^{-1}.$$

Therefore,

$$\|\overline{\mathbf{E}}_{l_n}(\beta)^{-1}\| = \lambda_{\max}[\overline{\mathbf{E}}_{l_n}(\beta)^{-1}] = \frac{1}{\lambda_{\min}[\overline{\mathbf{E}}_{l_n}(\beta)]} \leq m. \tag{40}$$

Using (39) and (40), we obtain:

$$\|\overline{\mathbf{E}}_{l_n}(\beta)^{-1}\mathbf{E}^*_{l_n-1}(\beta) - \mathbf{I}\| \leq \|\overline{\mathbf{E}}_{l_n}(\beta)^{-1}\| \cdot \|\mathbf{E}^*_{l_n-1}(\beta) - \overline{\mathbf{E}}_{l_n}(\beta)\| \to 0. \tag{41}$$

Comparing (38) and (41), we arrive at a contradiction.

The proof of (32) is very similar and is omitted. □

**Remark 3.10** Assume that $\overline{\mathbf{R}}^{(c)}_n(\beta) = \overline{\mathbf{R}}^{(c)}(\beta)$ for all $n$. In this case, condition (C) is satisfied by the sequence $\{\mathcal{R}^*_n(\beta)\}_n$ given in Example 2.4. By Theorem 3.9, the sequence (19) of estimating equations is an asymptotic quasi-score within the collection $\{\mathcal{H}_n\}_n$.

# 4 Asymptotic Existence and Strong Consistency

In this section, we fix a value $\beta_0 \in \mathcal{T}$, which we regard as the "true" (but unknown) value of the parameter $\beta$. Our aim is to give some sufficient conditions for the existence of a sequence of estimators $\hat{\beta}_n$, defined as solutions of (11), such that $\{\hat{\beta}_n\}_n$ converges to $\beta_0$ a.s. These conditions are slightly weaker than the conditions for the asymptotic optimality of the sequence $\{\mathbf{g}^*_n(\beta)\}_n$, encountered in Theorem 3.9. In particular, we may allow the matrix $\mathcal{R}^*_n(\beta)$ not to depend on $\beta$.

**Remark 4.1** In the present section, the a.s. statements refer to the probability $P_{\beta_0}$. Moreover, we employ the usual convention of omitting the argument $\beta_0$ in $\mathbf{g}_n^*(\beta_0), \mathbf{M}_n(\beta_0), \mathbf{H}_n(\beta_0)$, etc.

Recall that a sequence $\{\mathbf{s}_n\}_{n\geq 1}$ of $p$-dimensional random vectors with components $\mathbf{s}_n = (s_n^{(1)}, \ldots, s_n^{(p)})$, is called a *p-dimensional martingale* if

$$E(s_n^{(k)}|\mathbf{s}_1,\ldots,\mathbf{s}_{n-1}) = s_{n-1}^{(k)}, \quad \forall k \in \{1,\ldots,p\}, \ \forall n \geq 1.$$

Let $B_r = B_r(\beta_0) = \{\beta \in \mathcal{T}; \|\beta - \beta_0\| \leq r\}$ and $\partial B_r = \{\beta; \|\beta - \beta_0\| = r\}$.

We begin with a general result, which is similar to Theorem 7 of [27]. By way of comparison, we note that our result is formulated within a martingale context and that we do not require that the matrix $\mathcal{D}_n$ be positive definite.

**Theorem 4.2** *Let $\{\mathbf{q}_n(\beta), \beta \in \mathcal{T}\}_{n\geq 1}$ be a sequence of $p$-dimensional random functions, such that each $\mathbf{q}_n(\beta)$ is continuously differentiable with:*

$$\mathcal{D}_n(\beta) := -\frac{\partial \mathbf{q}_n(\beta)}{\partial \beta^T}, \quad \beta \in \mathcal{T}$$

*and $\{\mathbf{q}_n\}_{n\geq 1}$ is a $p$-dimensional martingale with mean zero and $\mathbf{M}_n = \mathrm{Cov}[\mathbf{q}_n]$. Let $\{\alpha_n\}_n$ be a sequence of constants such that, for some $C > 0$ and $N \geq 1$,*

$$\alpha_n \geq C\lambda_{\max}(\mathbf{M}_n), \quad \forall n \geq N. \tag{42}$$

*Assume that the following conditions hold:*

- (I) $\lambda_{\min}[\mathbf{M}_n] \to \infty$
- (S) *there exist some constants $\delta > 0, c_0 > 0$ such that, with probability $P_{\beta_0}$ equal to 1, there exists some random numbers $r_1 > 0, n_1 \geq 1$ for which*
  - (i) $|\lambda^T \mathcal{D}_n(\beta)\lambda| > 0$ for all $\lambda, \|\lambda\| = 1$, and for all $\beta \in B_{r_1}, n \geq n_1$;
  - (ii) $\lim_{r\to 0} \limsup_{n\to\infty} \alpha_n^{-1/2-\delta} \sup_{\beta \in B_r} \||\mathcal{D}_n(\beta) - \mathcal{D}_n\|| = 0$;
  - (iii) $|\lambda^T \mathcal{D}_n \lambda| \geq c_0 \alpha_n^{1/2+\delta}$ for all $\lambda, \|\lambda\| = 1$, and for all $n \geq n_1$.

*Then, there exists a sequence $\{\hat{\beta}_n\}_n \subset \mathcal{T}$ and a random number $n_0$ such that:*
  *(a) $P(\mathbf{q}_n(\hat{\beta}_n) = 0, \text{ for all } n \geq n_0) = 1$;*
  *(b) $\hat{\beta}_n \to \beta_0$ a.s.*

**Remark 4.3** Condition (S)(ii) says that, with probability $P_{\beta_0}$ equal to 1, the sequence $\{\alpha_n^{-1/2-\delta}\mathcal{D}_n(\beta)\}_{n\geq 1}$ is equicontinuous at $\beta_0$.

The proof of Theorem 4.2 combines analytic and stochastic techniques. On the analytic side, we have a result from topology, which provides an ingenious method for proving the existence of the solution of the one-to-one continuously differentiable function $\mathbf{q}_n(\beta)$. Breaking down condition $(S)$ into components has the advantage of allowing us to formulate some sufficient conditions for strong consistency of estimators, in terms of conditions on the matrices and functions that define our estimating equations (see Theorem 4.13). On the other hand, condition (I) enables us to apply the strong law of large numbers (SLLN) to the martingale $\{\mathbf{q}_n\}_{n \geq 1}$.

For the sake of completeness, we state some auxiliary results below.

**Lemma 4.4** *Let $T : \mathcal{T} \subset \mathbb{R}^p \to \mathbb{R}^p$ be a one-to-one continuously differentiable function and $\beta_0 \in \mathcal{T}$ such that $B_r(\beta_0) \subset \mathcal{T}$. If*

$$\|T(\beta_0)\| \leq \inf_{\beta \in \partial B_r} \|T(\beta) - T(\beta_0)\|,$$

*then there exists a (unique) $\hat{\beta} \in B_r$ such that $T(\hat{\beta}) = 0$.*

(This result is a consequence of Lemma A of [4].)

**Lemma 4.5 (Martingale SLLN)** *Let $\{\mathbf{s}_n\}_{n \geq 1}$ be a $p$-dimensional martingale with mean zero and covariance matrix $\mathbf{M}_n$. If $\lambda_{\min}(\mathbf{M}_n) \to \infty$, then*

$$\frac{\mathbf{s}_n}{[\lambda_{\max}(\mathbf{M}_n)]^{1/2+\delta}} \to 0 \quad a.s. \quad \forall \delta > 0.$$

(This result can be proved component-wise, using Theorem 3 of [15] with $p = 1$.)

**Proof of Theorem 4.2:** Let $\Omega_1$ be the event on which condition $(S)$ holds, with $P_{\beta_0}(\Omega_1) = 1$. For each $\omega \in \Omega_1$ and $n \geq n_1$, we consider the function

$$T_n(\beta) := \alpha_n^{-1/2-\delta} \mathbf{q}_n(\beta), \quad \beta \in B_{r_1}.$$

This function is continuously differentiable; it is also one-to-one, since $\mathcal{D}_n(\beta)$ is non-singular for all $\beta \in B_{r_1}$.

We prove that there exists an event $\Omega_0 \subset \Omega_1$ with $P_{\beta_0}(\Omega_0) = 1$, such that for every $\omega \in \Omega_0$ and for any $\varepsilon > 0$, there exist some random numbers $r_\varepsilon \in (0, \varepsilon), r_\varepsilon < r_1$ and $n_\varepsilon > n_1$ such that

$$\|T_n(\beta_0)\| \leq \inf_{\beta \in \partial B_{r_\varepsilon}} \|T_n(\beta) - T_n(\beta_0)\|, \quad \forall n \geq n_\varepsilon. \tag{43}$$

We claim that the conclusion of the theorem follows from here. To see this, note that by Lemma 4.4, on the event $\Omega_0$, for any $\varepsilon > 0$, there exists $\hat{\beta}_{n,\varepsilon} \in B_{r_\varepsilon}$ such that $T_n(\hat{\beta}_{n,\varepsilon}) = 0$ for all $n > n_\varepsilon$. Let $\varepsilon_0 > 0$ be fixed and denote $r_0 = r_{\varepsilon_0}$ and $n_0 = n_{\varepsilon_0}$. We define $\hat{\beta}_n := \hat{\beta}_{n,\varepsilon_0}$, for all $n \geq n_0$. Clearly, on the event $\Omega_0$, $T_n(\hat{\beta}_n) = 0$ for all $n \geq n_0$, i.e. (a) holds. If $\varepsilon > 0$ is arbitrary, then for any $n \geq n'_\varepsilon := \max\{n_\varepsilon, n_0\}$, both $\hat{\beta}_n$ and $\hat{\beta}_{n,\varepsilon}$ are zeros of the function $T_n$, in the ball $B_{r'_\varepsilon}$ with radius $r'_\varepsilon = \min\{r_\varepsilon, r_0\}$. Since $T_n$ is one-to-one, we conclude that $\hat{\beta}_{n,\varepsilon} = \hat{\beta}_n$ for all $n \geq n'_\varepsilon$. This argument shows that $\hat{\beta}_{n,\varepsilon}$ does not depend on $\varepsilon$, if $n$ is large enough. Finally, part (b) of the conclusion follows, since on the event $\Omega_0$, $\|\hat{\beta}_n - \beta_0\| \leq r_\varepsilon < \varepsilon$ for all $n \geq n_\varepsilon$.

We now turn to the proof of (43). We first treat the right-hand side. Let $\omega \in \Omega_1$ be fixed. By condition $(S)(ii)$, for any $\varepsilon > 0$, there exist some random numbers $r_\varepsilon \in (0, \varepsilon), r_\varepsilon < r_1$ and $N_\varepsilon > n_1$, such that for all $n > N_\varepsilon, \beta \in B_{r_\varepsilon}$,

$$\alpha_n^{-1/2-\delta} |||\mathcal{D}_n(\beta) - \mathcal{D}_n||| < \varepsilon.$$

and hence, for any vector $\lambda, \|\lambda\| = 1$, we have $\alpha_n^{-1/2-\delta}|\lambda^T[\mathcal{D}_n(\beta) - \mathcal{D}_n]\lambda| < \varepsilon$. In particular, it follows that for all $n > N_\varepsilon, \beta \in B_{r_\varepsilon}$,

$$\alpha_n^{-1/2-\delta}|\lambda^T \mathcal{D}_n(\beta)\lambda| > \alpha_n^{-1/2-\delta}|\lambda^T \mathcal{D}_n \lambda| - \varepsilon. \tag{44}$$

Let $n > N_\varepsilon$ and $\beta \in \partial B_{r_\varepsilon}$ be arbitrary. Using Taylor's formula, there exists $\bar{\beta}_n \in B_{r_\varepsilon}$ such that $\mathbf{q}_n(\beta) - \mathbf{q}_n = -\mathcal{D}_n(\bar{\beta}_n)(\beta - \beta_0)$. By letting $\lambda = (\beta - \beta_0)/r_\varepsilon$, we obtain:

$$\begin{aligned}\|\mathbf{q}_n(\beta) - \mathbf{q}_n\|^2 &= (\beta - \beta_0)^T \mathcal{D}_n(\bar{\beta}_n)^T \mathcal{D}_n(\bar{\beta}_n)(\beta - \beta_0) \\ &= \lambda^T \mathcal{D}_n(\bar{\beta}_n)^T \mathcal{D}_n(\bar{\beta}_n)\lambda r_\varepsilon^2 \geq [\lambda^T \mathcal{D}_n(\bar{\beta}_n)\lambda]^2 r_\varepsilon^2,\end{aligned}$$

where for the last inequality, we used the fact that $\lambda^T C^T C \lambda \geq (\lambda^T C \lambda)^2$ for any matrix $C$ and for any vector $\lambda$, with $\|\lambda\| = 1$ (Lemma 1, [27]). Taking the square-root, multiplying by $\alpha_n^{-1/2-\delta}$, and using (44) and $(S)(iii)$, we obtain:

$$\begin{aligned}\|T_n(\beta) - T_n(\beta_0)\| &\geq \alpha_n^{-1/2-\delta}|\lambda^T \mathcal{D}_n(\bar{\beta}_n)\lambda|r_\varepsilon \geq \{\alpha_n^{-1/2-\delta}|\lambda^T \mathcal{D}_n \lambda| - \varepsilon\} r_\varepsilon \\ &\geq (c_0 - \varepsilon)r_\varepsilon.\end{aligned}$$

Hence,

$$\inf_{\beta \in \partial B_{r_\varepsilon}} \|T_n(\beta) - T_n(\beta_0)\| \geq (c_0 - \varepsilon)r_\varepsilon, \quad \forall n > N_\varepsilon. \tag{45}$$

We now treat the left-hand of (43). By Lemma 4.5 (and using (42) and (I)),

$$\|T_n(\beta_0)\| = \alpha_n^{-1/2-\delta}\|\mathbf{q}_n\| \to 0, \quad \text{a.s.} \tag{46}$$

Denote by $\Omega_2$ the event where (46) holds. Let $\Omega_0 = \Omega_1 \cap \Omega_2$ and fix $\omega \in \Omega_0$. For any $\varepsilon > 0$, let $r_\varepsilon, N_\varepsilon$ be as above. By (46), there exists $n_\varepsilon > N_\varepsilon$ such that

$$\|T_n(\beta_0)\| \leq (c_0 - \varepsilon)r_\varepsilon, \quad \forall n > n_\varepsilon. \tag{47}$$

Relation (43) follows from (45) and (47). □

In what follows, we will apply the previous result to the case of the estimating function $\mathbf{g}_n^*(\beta)$. As in [27], we have:

$$\mathcal{D}_n^*(\beta) := -\frac{\partial \mathbf{g}_n^*(\beta)}{\partial \beta^T} = \mathbf{H}_n(\beta) - \mathbf{B}_n(\beta) - \mathcal{E}_n(\beta), \quad \beta \in \mathcal{T}$$

(see Appendix A for the exact formulas of $\mathbf{H}_n(\beta), \mathbf{B}_n(\beta), \mathcal{E}_n(\beta)$).

We define the following constants:

$$\begin{aligned}\gamma_n^{(0),\text{indep}} &= \max_{i \leq n, j \leq m} \mathbf{x}_{ij}^T (\mathbf{H}_n^{\text{indep}})^{-1} \mathbf{x}_{ij} \\ a_n &= \lambda_{\max}(\mathbf{H}_n^{\text{indep}}) \gamma_n^{(0),\text{indep}}.\end{aligned}$$

Suppose that $\mu$ is three times continuously differentiable. For any $r > 0$ and $n \geq 1$, we let

$$k_n^{[2]}(r) = \sup_{\beta \in B_r} \max_{i \leq n, j \leq m} \left| \frac{\mu''(\mathbf{x}_{ij}^T \beta)}{\mu'(\mathbf{x}_{ij}^T \beta)} \right|, \quad k_n^{[3]}(r) = \sup_{\beta \in B_r} \max_{i \leq n, j \leq m} \left| \frac{\mu'''(\mathbf{x}_{ij}^T \beta)}{\mu'(\mathbf{x}_{ij}^T \beta)} \right|$$

$$\eta_n(r) = \sup_{\beta, \beta' \in B_r} \max_{i \leq n, j \leq m} \left| \left[ \frac{\mu'(\mathbf{x}_{ij}^T \beta')}{\mu'(\mathbf{x}_{ij}^T \beta)} \right]^{1/2} - 1 \right|$$

$$\pi_n(r) = \sup_{\beta \in B_r} \max_{i \leq n} \lambda_{\max}[(\mathcal{R}_{i-1}^*)^{1/2} \mathcal{R}_{i-1}^*(\beta)^{-1} (\mathcal{R}_{i-1}^*)^{1/2}]$$

$$\rho_n(r) = \sup_{\beta \in B_r} \max_{i \leq n} \lambda_{\max}[(\mathcal{R}_{i-1}^*)^{1/2} \mathcal{R}_{i-1}^*(\beta)^{-1} (\mathcal{R}_{i-1}^*)^{1/2} - \mathbf{I}]$$

$$q_n(r) = \sup_{\beta \in B_r} \max_{i \leq n, l \leq p} \lambda_{\max} \left[ \frac{\partial}{\partial \beta_l} \mathcal{R}_{i-1}^*(\beta) \right]$$

As in [27], we introduce the following assumption:

$(AH)$     there exists $C > 0, r_0 > 0$ such that $k_n^{[l]}(r_0) \leq C$ for all $n \geq 1, l = 2, 3$.

Note that $(AH)$ holds if the covariates are bounded.

We introduce a new assumption:

$$(K) \quad \lim_{r \to 0} \limsup_{n \to \infty} r a_n^{1/2} = 0.$$

We have the following result, whose proof is given in Appendix B.

**Lemma 4.6** *Under $(AH)$ and $(K)$, there exist $r_1 > 0$ and $n_1 \geq 1$ such that*

$$\eta_n(r) \leq C r a_n^{1/2}, \quad \text{for all } r \in (0, r_1), n \geq n_1.$$

*In particular, under $(AH)$ and $(K)$, $\lim_{r \to 0} \limsup_{n \to \infty} \eta_n(r) = 0$.*

We consider the following condition on the sequence $\{\mathcal{R}_n^*(\beta)\}_n$:

$(R')$     there exists $C > 0$ such that $\lambda_{\min}(\mathcal{R}_n^*) \geq C$ for all $n \geq 1$ a.s.

Clearly, condition $(R')$ is weaker than condition $(R)$ (encountered in Theorem 3.9). The following fact is an immediate consequence of condition $(R')$:

$$\lambda_{\max}[\mathcal{R}_{i-1}^*(\beta)^{-1}] \leq C \pi_n(r), \quad \forall \beta \in B_r, \forall i \leq n, \text{ a.s.} \tag{48}$$

**Remark 4.7** Suppose that with probability $P_{\beta_0}$ equal to 1, the sequence $\{\mathcal{R}_n^*(\beta)\}_{n \geq 1}$ is equicontinuous at $\beta_0$, i.e.

$$\lim_{r \to 0} \limsup_{n \to \infty} \sup_{\beta \in B_r} \|\mathcal{R}_n^*(\beta) - \mathcal{R}_n^*\| = 0 \quad \text{a.s.} \tag{49}$$

If the sequence $\{\mathcal{R}_n^*(\beta)\}_n$ satisfies $(R')$, then $\lim_{r \to 0} \limsup_{n \to \infty} \rho_n(r) = 0$ a.s. and $\lim_{r \to 0} \limsup_{n \to \infty} \pi_n(r) = 1$ a.s.

**Example 4.8** Assume that (15) holds. Let $\{\mathcal{R}_n^*(\beta)\}_n$ be the sequence introduced in Example 2.4. Using the same argument as in the proof of Proposition 2, [2], one can show that if:

(i) $\lim_{r \to 0} \limsup_{n \to \infty} \eta_n(r) = 0$,
(ii) $\lambda_{\max}(\mathbf{H}_n^{\text{indep}}) \leq Cn$ for all $n \geq 1$,
(iii) $E\|\mathbf{A}_i^{-1/2}\varepsilon_i\|^{2+\delta} \leq C$ for all $i \geq 1$, for some $\delta > 0$

then (49) holds.

Let $\{\alpha_n\}_n$ be a non-decreasing sequence of constants with $\lim_n \alpha_n = \infty$, and

$$\delta_n = \alpha_n^{-1/2-\delta} \lambda_{\max}(\mathbf{H}_n^{\text{indep}}).$$

The following three lemmas examine the asymptotic behavior of the three terms of $\mathcal{D}_n^*(\beta)$. Their respective proofs are given in Appendices C, D and E (see also [1]).

**Lemma 4.9** *Suppose that* $(AH)$ *and* $(K)$ *hold. Let* $\{\mathcal{R}_n^*(\beta)\}_n$ *be a sequence of random matrices which satisfy* $(A)$, $(B)$ *and* $(R')$. *If*

$(C_1')$ $\quad \lim_{r \to 0} \limsup_{n \to \infty} \delta_n \pi_n(r) \eta_n(r) = 0$ a.s.

$(C_2)$ $\quad \lim_{r \to 0} \limsup_{n \to \infty} \delta_n \rho_n(r) = 0$ a.s.

then

$$\lim_{r \to 0} \limsup_{n \to \infty} \alpha_n^{-1/2-\delta} \sup_{\beta \in B_r} \|\|\mathbf{H}_n(\beta) - \mathbf{H}_n\|\| = 0 \quad a.s.$$

**Lemma 4.10** *Suppose that* $(AH)$ *and* $(K)$ *hold. Let* $\{\mathcal{R}_n^*(\beta)\}_n$ *be a sequence of random matrices which satisfy* $(A)$, $(B)$ *and* $(R')$. *If*

$(C_1)$ $\quad \lim_{r \to 0} \limsup_{n \to \infty} r\delta_n \pi_n(r) a_n^{1/2} = 0$ a.s.

$(C_3)$ $\quad \lim_{r \to 0} \limsup_{n \to \infty} r\delta_n \pi_n^2(r) q_n(r) = 0$ a.s.,

then

$$\lim_{r \to 0} \limsup_{n \to \infty} \alpha_n^{-1/2-\delta} \sup_{\beta \in B_r} \|\|\mathbf{B}_n(\beta)\|\| = 0 \quad a.s.$$

**Lemma 4.11** *Suppose that* $(AH)$ *and* $(K)$ *hold. Let* $\{\mathcal{R}_n^*(\beta)\}_n$ *be a sequence of random matrices which satisfy* $(A)$, $(B)$ *and* $(R')$. *If*

$(C_4)$ $\quad \limsup_{n \to \infty} nE[\pi_n^2(r)]\tilde{a}_n \lambda_{\max}(\mathbf{H}_n^{\text{indep}}) < \infty$, where $\tilde{a}_n = \max\{a_n, a_n^2\}$

$(C_5)$ $\quad \limsup_{n \to \infty} nE[\pi_n^4(r) q_n^2(r)] \lambda_{\max}(\mathbf{H}_n^{\text{indep}}) < \infty$ for all $r > 0$,

then

$$\lim_{n \to \infty} \alpha_n^{-1/2-\delta} \sup_{\beta \in B_r} \|\|\mathcal{E}_n(\beta)\|\| = 0 \quad a.s.$$

**Remark 4.12** 1. In the case of the linear regression model, $\mu(x) = x$ for all $x$. Hence $\eta_n(r) = 0$ for all $r > 0, n \geq 1$, and $(C_1')$ is automatically satisfied. In this case $\mathbf{A}_i = \mathbf{I}$ and $\mathbf{H}_n^{\text{indep}} = \sum_{i=1}^n \mathbf{X}_i^T \mathbf{X}_i$.
2. Lemma 4.6 shows that condition $(C_1)$ is stronger that $(C_1')$.
3. If $\mathcal{R}_n^*(\beta)$ does not depend on $\beta$, then $\rho_n(r) = q_n(r) = 0$ for all $r > 0, n \geq 1$; hence, $(C_2)$, $(C_3)$ and $(C_5)$ are satisfied. In particular, this is the case of Examples 2.1, 2.2, and 2.3.

As in [2], we introduce the following assumption:

$(H')$    there exists a constant $C > 0$ such that $\lambda_{\min}(\overline{\mathbf{R}}_n^{(c)}) \geq C, \forall n \geq 1$, a.s.

Here is the main result of this section.

**Theorem 4.13** *Suppose that $(H')$, $(AH)$ and $(K)$ hold. Let $\{\mathcal{R}_n^*(\beta)\}_n$ be a sequence of random matrices which satisfy $(A)$, $(B)$, $(R')$ and*

$(E)$    *there exists a constant $C > 0$ such that $\lambda_{\max}(\mathcal{R}_n^*) \leq C, \forall n \geq 1$, a.s.*

*Suppose that conditions $(C_1)$-$(C_5)$ are satisfied with $\alpha_n = \lambda_{\max}(\mathbf{H}_n^{\text{indep}})$. If*
*(i) $\lambda_{\min}(\mathbf{H}_n^{\text{indep}}) \to \infty$*
*(ii) there exist an integer $N \geq 1$ and some constants $\delta > 0, c_0 > 0$ such that*

$$\lambda_{\min}(\mathbf{H}_n^{\text{indep}}) \geq c_0 [\lambda_{\max}(\mathbf{H}_n^{\text{indep}})]^{1/2+\delta}, \quad \forall n \geq N,$$

*then there exists a sequence $\{\hat{\beta}_n\}_n \subset \mathcal{T}$ and a random number $n_0$ such that:*
*(a) $P(\mathbf{g}_n^*(\hat{\beta}_n) = 0, \text{ for all } n \geq n_0) = 1$;*
*(b) $\hat{\beta}_n \to \beta_0$ a.s.*

**Remark 4.14** Hypothesis $(i)$ and $(ii)$ of Theorem 4.13 are indeed very mild. To see this, consider the following stronger form of hypothesis $(ii)$:

$(ii)'$ $\lambda_{\min}[\mathbf{H}_n^{\text{indep}}(\beta)] \geq c_0 [\lambda_{\max}(\mathbf{H}_n^{\text{indep}})]^{1/2+\delta}, \quad \beta \in \mathcal{T}, \forall n \geq N.$

Using the approach of [6], one can prove that under $(i)$ and $(ii)'$, there exists a sequence of strongly consistent estimators, defined as roots of the "working independence" equation $\mathbf{g}_n^{\text{indep}}(\beta) = 0$ (see Remark 3.2).

**Remark 4.15** (Discussion of hypothesis $(i)$ and $(ii)$ of Theorem 4.13) In the case of the usual regression models, conditions $(i)$ and $(ii)$ of Theorem 4.13 can be simplified into conditions which speak only about the asymptotic behavior of the covariates $(\mathbf{x}_{ij})_{i \geq 1}$ for $j = 1, \ldots, m$.
   1. In the case of the linear regression model, $\mu(x) = x$ and $\mathbf{H}_n^{\text{indep}} = \sum_{i=1}^n \sum_{j=1}^m \mathbf{x}_{ij} \mathbf{x}_{ij}^T$. Hypothesis $(i)$ holds if $\sum_{i \geq 1} \lambda_{\min}(\mathbf{x}_{ij} \mathbf{x}_{ij}^T) = \infty$ for some $j$, whereas $(ii)$ holds if there exist an integer $N \geq 1$ and some constants $\delta > 0, c_0 > 0$ such that

$$\sum_{i=1}^n \sum_{j=1}^m \lambda_{\min}(\mathbf{x}_{ij} \mathbf{x}_{ij}^T) \geq c_0 \left[ \sum_{i=1}^n \sum_{j=1}^m \lambda_{\max}(\mathbf{x}_{ij} \mathbf{x}_{ij}^T) \right]^{1/2+\delta}. \tag{50}$$

If the covariates $(\mathbf{x}_{ij})_{i\geq 1}$ are bounded, then (50) holds if there exists some $j = 1, \ldots, m$ such that $\sum_{i=1}^{n} \lambda_{\min}(\mathbf{x}_{ij}\mathbf{x}_{ij}^T) \geq c_0 n^{1/2+\delta}$ for all $n \geq N$.

2. In the case of the logistic regression model, $\mu(x) = \exp(x)/[1 + \exp(x)]$ and $\mathbf{H}_n^{\text{indep}} = \sum_{i=1}^{n} \sum_{j=1}^{m} \exp\{\mathbf{x}_{ij}^T \beta_0\}/(1 + \exp\{\mathbf{x}_{ij}^T \beta_0\})^2 \mathbf{x}_{ij}\mathbf{x}_{ij}^T$. Assuming that the parameter set $\mathcal{T}$ is bounded by a constant $C > 0$, and letting $a_{ij} = \exp\{-C\|\mathbf{x}_{ij}\|\}/(1 + \exp\{C\|\mathbf{x}_{ij}^T\|\})$, we see that

$$\frac{1}{2}\sum_{i=1}^{n}\sum_{j=1}^{m} a_{ij}\mathbf{x}_{ij}\mathbf{x}_{ij}^T \leq \mathbf{H}_n^{\text{indep}} \leq \sum_{i=1}^{n}\sum_{j=1}^{m} \mathbf{x}_{ij}\mathbf{x}_{ij}^T.$$

Hypothesis $(i)$ holds if $\sum_{i\geq 1} a_{ij}\lambda_{\min}(\mathbf{x}_{ij}\mathbf{x}_{ij}^T) = \infty$ for some $j$, whereas $(ii)$ holds if there exist an integer $N \geq 1$ and some constants $\delta > 0, c_0 > 0$ such that

$$\sum_{i=1}^{n}\sum_{j=1}^{m} a_{ij}\lambda_{\min}(\mathbf{x}_{ij}\mathbf{x}_{ij}^T) \geq c_0 \left[\sum_{i=1}^{n}\sum_{j=1}^{m} \lambda_{\max}(\mathbf{x}_{ij}\mathbf{x}_{ij}^T)\right]^{1/2+\delta}, \quad \forall n \geq N. \quad (51)$$

If the covariates $(\mathbf{x}_{ij})_{i\geq 1}$ are bounded, then (51) holds if there exists some $j = 1, \ldots, m$ such that $\sum_{i=1}^{n} a_{ij}\lambda_{\min}(\mathbf{x}_{ij}\mathbf{x}_{ij}^T) \geq c_0 n^{1/2+\delta}$ for all $n \geq N$. (See also [7] for a related analysis, in the case $m = 1$.)

3. In the case of the Poisson regression model, $\mu(x) = \exp(x)$ and $\mathbf{H}_n^{\text{indep}} = \sum_{i=1}^{n}\sum_{j=1}^{m} \exp\{\mathbf{x}_{ij}^T\beta_0\}\mathbf{x}_{ij}\mathbf{x}_{ij}^T$. Assume that the parameter set $\mathcal{T}$ is bounded by a constant $C > 0$, and let $b_{ij} = \exp\{C\|\mathbf{x}_{ij}^T\|\})$. Then

$$\sum_{i=1}^{n}\sum_{j=1}^{m} \frac{1}{b_{ij}}\mathbf{x}_{ij}\mathbf{x}_{ij}^T \leq \mathbf{H}_n^{\text{indep}} \leq \sum_{i=1}^{n}\sum_{j=1}^{m} b_{ij}\mathbf{x}_{ij}\mathbf{x}_{ij}^T.$$

Hypothesis $(i)$ holds if $\sum_{i\geq 1} \lambda_{\min}(\mathbf{x}_{ij}\mathbf{x}_{ij}^T)/b_{ij} = \infty$ for some $j$, whereas $(ii)$ holds if there exist an integer $N \geq 1$ and some constants $\delta > 0, c_0 > 0$ such that

$$\sum_{i=1}^{n}\sum_{j=1}^{m} \frac{1}{b_{ij}}\lambda_{\min}(\mathbf{x}_{ij}\mathbf{x}_{ij}^T) \geq c_0 \left[\sum_{i=1}^{n}\sum_{j=1}^{m} b_{ij}\lambda_{\max}(\mathbf{x}_{ij}\mathbf{x}_{ij}^T)\right]^{1/2+\delta}, \quad \forall n \geq N. \quad (52)$$

If the covariates $(\mathbf{x}_{ij})_{i\geq 1}$ are bounded, then (52) holds if there exists some $j = 1, \ldots, m$ such that $\sum_{i=1}^{n} \lambda_{\min}(\mathbf{x}_{ij}\mathbf{x}_{ij}^T) \geq c_0 n^{1/2+\delta}$ for all $n \geq N$.

**Proof of Theorem 4.13:** We will apply Theorem 4.2 to the function $\mathbf{q}_n(\beta) = \mathbf{g}_n^*(\beta)$, by taking $\alpha_n = \lambda_{\max}(\mathbf{H}_n^{\text{indep}})$. Due to $(R')$, $E|\mathbf{g}_n^*| < \infty$, for any $n$. Hence, $\{\mathbf{g}_n^*\}_n$ is a martingale. Recall that $\mathbf{M}_n^* = \sum_{i=1}^{n} \mathbf{X}_i \mathbf{A}_i^{1/2} \overline{\mathbf{E}}_{i-1}^* \mathbf{A}_i^{1/2} \mathbf{X}_i$, where $\overline{\mathbf{E}}_{i-1}^* = E[(\mathcal{R}_{i-1}^*)^{-1}\overline{\mathbf{R}}_i^{(c)}(\mathcal{R}_{i-1}^*)^{-1}]$ (see (20)).

We first prove that (42) holds. Using $(R')$ and the fact that $\overline{\mathbf{R}}_i^{(c)} \leq m\mathbf{I}$ for all $i \geq 1$ a.s., it follows that $(\mathcal{R}_{i-1}^*)^{-1}\overline{\mathbf{R}}_i^{(c)}(\mathcal{R}_{i-1}^*)^{-1} \leq C\mathbf{I}$ for all $i \geq 1$, a.s. Hence $\overline{\mathbf{E}}_i^* \leq C\mathbf{I}$ for all $i \geq 1$ and

$$\mathbf{M}_n^* \leq C\mathbf{H}_n^{\text{indep}}, \quad \forall n \geq 1.$$

Hence $\lambda_{\max}(\mathbf{M}_n^*) \leq C\lambda_{\max}(\mathbf{H}_n^{\text{indep}})$, i.e. $\alpha_n$ satisfies relation (42).

We now prove that $(I)$ holds. For any $p \times 1$ vector $\lambda$ with $\|\lambda\| = 1$, we have:

$$\sum_{i=1}^n \lambda^T \mathbf{X}_i^T \mathbf{A}_i^{1/2} (\mathcal{R}_{i-1}^*)^{-1} \overline{\mathbf{R}}_i^{(c)} (\mathcal{R}_{i-1}^*)^{-1} \mathbf{A}_i^{1/2} \mathbf{X}_i \lambda \geq \min_{i \leq n} \lambda_{\min}(\overline{\mathbf{R}}_i^{(c)}) \cdot$$

$$\min_{i \leq n} \lambda_{\min}[(\mathcal{R}_{i-1}^*)^{-2}] \cdot \lambda^T \mathbf{H}_n^{\text{indep}} \lambda \geq C_0 \lambda^T \mathbf{H}_n^{\text{indep}} \lambda,$$

using $(H')$ and $(E)$. Taking the expectation (with respect to $P_{\beta_0}$), we conclude that $\lambda^T \mathbf{M}_n^* \lambda \geq C_0 \lambda^T \mathbf{H}_n^{\text{indep}} \lambda$, i.e.

$$\mathbf{M}_n^* \geq C_0 \mathbf{H}_n^{\text{indep}}, \quad \forall n \geq 1.$$

The fact that $(I)$ holds follows from our hypothesis (i).

We now prove that $(S)$ holds. Part $(ii)$ follows directly from Lemmas 4.9, 4.10 and 4.11. To prove that parts $(i)$ and $(iii)$ hold, note that condition $(E)$ and our hypothesis $(ii)$ imply that for any $p \times 1$ vector $\lambda$ with $\|\lambda\| = 1$:

$$\lambda^T \mathbf{H}_n \lambda = \sum_{i=1}^n \lambda^T \mathbf{X}_i^T \mathbf{A}_i^{1/2} (\mathcal{R}_{i-1}^*)^{-1} \mathbf{A}_i^{1/2} \mathbf{X}_i \lambda \geq \min_{i \leq n} \lambda_{\min}[(\mathcal{R}_{i-1}^*)^{-1}] \cdot \lambda^T \mathbf{H}_n^{\text{indep}} \lambda$$

$$\geq C \lambda^T \mathbf{H}_n^{\text{indep}} \lambda \geq C \alpha_n^{1/2+\delta}, \; \forall n \geq N.$$

From Lemmas 4.9, 4.10, 4.11, it follows that, with probability 1, there exist some random numbers $r_1 > 0, n_1 \geq 1$ such that

$$\alpha_n^{-1/2-\delta} |\lambda^T [\mathbf{H}_n(\beta) - \mathbf{H}_n]\lambda + \lambda^T [\mathbf{B}_n(\beta) + \mathcal{E}_n(\beta)]\lambda| < C/2, \quad \forall \beta \in B_{r_1}, n \geq n_1.$$

Recalling that $\mathcal{D}_n(\beta) = \mathbf{H}_n(\beta) - \mathbf{B}_n(\beta) - \mathcal{E}_n(\beta)$, we conclude that:

$$\begin{aligned}|\lambda^T \mathcal{D}_n(\beta)\lambda| &\geq |\lambda^T \mathbf{H}_n \lambda| - |\lambda^T[\mathbf{H}_n(\beta) - \mathbf{H}_n]\lambda + \lambda^T[\mathbf{B}_n(\beta) + \mathcal{E}_n(\beta)]\lambda| \\ &\geq C\alpha_n^{1/2+\delta} - C\alpha_n^{1/2+\delta}/2 > 0, \quad \text{for all } \beta \in B_{r_1}, n \geq n_1,\end{aligned}$$

This concludes the proof of parts $(i)$ and $(iii)$ of $(S)$. $\square$

**Remark 4.16** The proof of Theorem 4.13 can be adapted to apply to the sequence $\{\mathbf{g}_n^{\text{indep}}(\beta)\}_{n \geq 1}$ (in fact, only the proof of Lemma 4.9 needs to be adapted, since $\partial \mathbf{g}_n^{\text{indep}}(\beta)/\partial \beta^T = \mathbf{H}_n^{\text{indep}}(\beta)$). More precisely, assume that $(H')$ holds and

$$(K') \quad \lim_{r \to 0} \limsup_{n \to \infty} \eta_n(r) = 0.$$

Under conditions $(i)$ and $(ii)$ of Theorem 4.13, one can prove that there exist a sequence $\{\tilde{\beta}_n\}_n$ and a random integer $n_0$ such that

$$P(\mathbf{g}_n^{\text{indep}}(\tilde{\beta}_n) = 0, \forall n \geq n_0) = 1 \quad \text{and} \quad \tilde{\beta}_n \to \beta_0 \text{ a.s}$$

Note that $(K')$ holds if the covariates are bounded, or we have a linear regression model. Our set-up covers a more general situation than Theorem 2 of [6]; in our case, neither the joint distribution of the data nor the covariance matrices are known.

**Remark 4.17** Suppose that
$$\alpha_n = \lambda_{\max}(\mathbf{H}_n^{\mathrm{indep}}).$$

Then $\delta_n = [\lambda_{\max}(\mathbf{H}_n^{\mathrm{indep}})]^{1/2-\delta}$ and conditions $(C_1)$-$(C_3)$ become:

$(C_1^*)$ $\quad \lim_{r\to 0} \limsup_{n\to\infty} r\pi_n(r)(\gamma_n^{(0),\mathrm{indep}})^{1/2}[\lambda_{\max}(\mathbf{H}_n^{\mathrm{indep}})]^{1-\delta} = 0$ a.s

$(C_2^*)$ $\quad \lim_{r\to 0} \limsup_{n\to\infty} \rho_n(r)[\lambda_{\max}(\mathbf{H}_n^{\mathrm{indep}})]^{1/2-\delta} = 0$ a.s

$(C_3^*)$ $\quad \lim_{r\to 0} \limsup_{n\to\infty} r\pi_n^2(r)q_n(r)[\lambda_{\max}(\mathbf{H}_n^{\mathrm{indep}})]^{1/2-\delta} = 0$ a.s.

**Remark 4.18** (Discussion of the assumptions on $\{\mathcal{R}_n^*(\beta)\}_n$ in Theorem 4.13) In the case of Examples 2.2-2.4, the assumptions imposed in Theorem 4.13 on the sequence $\{\mathcal{R}_n^*(\beta)\}_n$ can be summarized as follows:

| Example | Assumptions |
|---|---|
| 2.2 (GEE) | $(C_1^*), (C_4), (R'), (E)$ |
| 2.3 (PLE) | $(C_1^*), (C_4), (R'), (E)$ |
| 2.4 (AQS) | $(C_1^*), (C_2^*), (C_3^*), (C_4), (C_5), (R'), (E)$ |

We note that assumption $(E)$ is implied by:

$(C')$ $\quad \mathcal{R}_{n-1}^* - \overline{\mathbf{R}}_n^{(c)} \to 0$ (element-wise) a.s.,

which is satisfied by the sequences $\{\mathcal{R}_n^*(\beta)\}_n$ given in Examples 2.3 and 2.4.

**Remark 4.19** (Weak consistency and asymptotic normality) Let $\beta_0$ be the true value of the parameter, $\mathcal{D}_n^*(\beta) = \partial \mathbf{g}_n^*(\beta)/\partial \beta^T$, $\mathbf{M}_n^* = \mathbf{M}_n^*(\beta_0)$, and $\mathbf{H}_n^* = \mathbf{H}_n^*(\beta_0)$. Using a methodology similar to [2] and [27], if we let

$$\tau_n^* = m \max_{i\le n} \lambda_{\max}[(\mathcal{R}_{i-1}^*)^{-1}], \quad c_n^* = \lambda_{\max}[(\mathbf{M}_n^*)^{-1}\mathbf{H}_n^*],$$

and $B_n^*(r) = \{\beta \in \mathcal{T}; \|(\mathbf{H}_n^*)^{1/2}(\hat{\beta}_n - \beta_0)\| \le (\tau_n^*)^{1/2}r\}$, then under the condition that $\{c_n^*\tau_n^*\}_n$ is bounded, and

$(CC^*)$ $\quad \sup_{\beta \in B_n^*(r)} \|(\mathbf{H}_n^*)^{-1/2}\mathcal{D}_n^*(\beta)(\mathbf{H}_n^*)^{-1/2} - \mathbf{I}\| \xrightarrow{P_{\beta_0}} 0 \quad \forall r > 0,$

one can prove that there exists a sequence $\{\hat{\beta}_n\}_n$ of weakly consistent estimators of $\beta_0$, such that $P_{\beta_0}(\mathbf{g}_n^*(\hat{\beta}_n) = 0) \to 1$ and

$$(\mathbf{M}_n^*)^{-1/2}\mathbf{g}_n^*(\beta_0) = (\mathbf{M}_n^*)^{-1/2}\mathbf{H}_n^*(\hat{\beta}_n - \beta_0) + o_{P_{\beta_0}}(1).$$

By an invocation of a martingale central limit theorem, under the appropriate conditions, one can conclude that $(\mathbf{M}_n^*)^{-1/2}\mathbf{g}_n^*(\beta_0) \xrightarrow{d} N(0, \mathbf{I})$, and therefore

$$(\mathbf{M}_n^*)^{-1/2}\mathbf{H}_n^*(\hat{\beta}_n - \beta_0) \xrightarrow{d} N(0, \mathbf{I}).$$

In view of (31) and (32), the matrices $\mathbf{M}_n^*$ and $\mathbf{H}_n^*$ are asymptotically the same, and hence
$$(\mathbf{H}_n^*)^{1/2}(\hat{\beta}_n - \beta_0) \xrightarrow{d} N(0, \mathbf{I}).$$

The matrix $\mathbf{H}_n^*$ depends on the unknown parameter $\beta_0$; it also depends on the matrix $\mathcal{R}_{i-1}^*(\beta_0)$ through the value $\mathbf{E}_{i-1}^*(\beta_0) = E_{\beta_0}[\mathcal{R}_{i-1}^*(\beta_0)^{-1}]$, which cannot be calculated from the data. As suggested by Remark 8, [27], in practice, one may approximate the matrix $\mathbf{H}_n^*$ by the matrix
$$\widehat{\mathbf{H}}_n := \sum_{i=1}^n \mathbf{X}_i^T \mathbf{A}_i(\hat{\beta}_n)^{1/2} \mathcal{R}_{i-1}^*(\hat{\beta}_n)^{-1} \mathbf{A}_i(\hat{\beta}_n)^{1/2} \mathbf{X}_i,$$
and obtain a confidence interval for $\beta_0$. (We do not discuss here the theoretical issues related to this practical implementation.)

**Remark 4.20** (The linear regression model) We consider separately the 4 estimating equations introduced in Examples 2.1-2.4, in the case of the linear regression model (i.e. $\mu(x) = x$).

1. (Working independence) The equation introduced in Example 2.1 has the solution
$$\hat{\beta}_n^{\text{indep}} = \left( \sum_{i=1}^n \mathbf{X}_i^T \mathbf{X}_i \right)^{-1} \left( \sum_{i=1}^n \mathbf{X}_i^T \mathbf{y}_i \right)$$
and asymptotic covariance matrix $\mathbf{H}_n^{\text{indep}} = \sum_{i=1}^n \mathbf{X}_i^T \mathbf{X}_i$.

2. (GEE) The equation introduced in Example 2.2 has the solution
$$\hat{\beta}_n^{\text{GEE}} = \left( \sum_{i=1}^n \mathbf{X}_i^T \mathbf{R}_i(\alpha)^{-1} \mathbf{X}_i \right)^{-1} \left( \sum_{i=1}^n \mathbf{X}_i^T \mathbf{R}_i(\alpha)^{-1} \mathbf{y}_i \right)$$
and asymptotic covariance matrix $\widehat{\mathbf{H}}_n^{\text{GEE}} = \sum_{i=1}^n \mathbf{X}_i^T \mathbf{R}_i(\alpha)^{-1} \mathbf{X}_i$. (The matrix $\mathbf{R}_i(\alpha)$ is supposed to be known.)

3. (Pseudo-likelihood equation) Let $\widetilde{\mathcal{R}}_0 = \widetilde{\mathcal{R}}_1 = \mathbf{I}$ and
$$\widetilde{\mathcal{R}}_k := \frac{1}{k} \sum_{i=1}^k (\mathbf{y}_i - \mathbf{X}_i \hat{\beta}_n^{\text{indep}})(\mathbf{y}_i - \mathbf{X}_i \hat{\beta}_n^{\text{indep}})^T, \quad k = 2, \ldots, n.$$

The equation introduced in Example 2.3 has the solution
$$\hat{\beta}_n^{\text{PLE}} = \left( \sum_{i=1}^n \mathbf{X}_i^T \widetilde{\mathcal{R}}_{i-1}^{-1} \mathbf{X}_i \right)^{-1} \left( \sum_{i=1}^n \mathbf{X}_i^T \widetilde{\mathcal{R}}_{i-1}^{-1} \mathbf{y}_i \right)$$
and asymptotic covariance matrix $\widehat{\mathbf{H}}_n^{\text{PLE}} = \sum_{i=1}^n \mathbf{X}_i^T \widetilde{\mathcal{R}}_{i-1}^{-1} \mathbf{X}_i$.

4. (Asymptotic Quasi-Score) Let $\mathcal{R}_0^* = \mathcal{R}_1^* = \mathbf{I}$ and
$$\mathcal{R}_k^*(\beta) := \frac{1}{k} \sum_{i=1}^k (\mathbf{y}_i - \mathbf{X}_i \beta)(\mathbf{y}_i - \mathbf{X}_i \beta)^T, \quad k = 2, \ldots, n.$$

The equation introduced in Example 2.4 is a polynomial of degree $(2m)^{n-2}+1$ in $\beta$. We select $\hat{\beta}_n^{AQS}$ to be the root of this polynomial, which is closest to $\hat{\beta}_n^{\text{indep}}$. This root cannot be written in closed form. The asymptotic covariance matrix of $\hat{\beta}_n^{\text{AQS}}$ is $\hat{\mathbf{H}}_n^{\text{AQS}} = \sum_{i=1}^n \mathbf{X}_i^T \mathcal{R}_{i-1}^*(\hat{\beta}_n^{\text{AQS}})^{-1}\mathbf{X}_i$.

The fact that the sequence $\{\mathbf{g}_n^*(\beta)\}_{n\geq 1}$ (given by Example 2.4) is an AQS tells us that, if $n$ is sufficiently large, then for any $k = 1,\ldots,p$, the asymptotic variance of $\hat{\beta}_{k,n}^{\text{AQS}}$ (i.e. the $(k,k)$-element of the matrix $\hat{\mathbf{H}}_n^{\text{AQS}}$) is smaller than the asymptotic variance of each of $\hat{\beta}_{k,n}^{\text{indep}}$, $\hat{\beta}_{k,n}^{\text{GEE}}$ and $\hat{\beta}_{k,n}^{\text{PLE}}$. (Here, we used the notation $\hat{\beta}_n = (\hat{\beta}_{1,n}, \ldots, \hat{\beta}_{p,n})^T$ in each of the 4 cases.)

# A    Formulas for the terms of $\mathcal{D}_n^*(\beta)$

We write $\mathbf{g}_n^*(\beta) = \mathbf{g}_{n,1}^*(\beta) + \mathbf{g}_{n,2}^*(\beta)$, where

$$\mathbf{g}_{n,1}^*(\beta) = \sum_{i=1}^n \mathbf{X}_i^T \mathbf{A}_i(\beta)^{1/2} \mathcal{R}_{i-1}^*(\beta)^{-1} \mathbf{A}_i(\beta)^{-1/2} [\mu_i - \mu_i(\beta)]$$

$$\mathbf{g}_{n,2}^*(\beta) = \sum_{i=1}^n \mathbf{X}_i^T \mathbf{A}_i(\beta)^{1/2} \mathcal{R}_{i-1}^*(\beta)^{-1} \mathbf{A}_i(\beta)^{-1/2} \varepsilon_i.$$

Note that

$$\frac{\partial \mathbf{g}_{n,1}^*(\beta)}{\partial \beta^T} = \mathbf{B}_n^{[1]}(\beta) + \mathbf{B}_n^{[2]}(\beta) + \mathbf{B}_n^{[3]}(\beta) - \mathbf{H}_n(\beta) := \mathbf{B}_n(\beta) - \mathbf{H}_n(\beta)$$

$$\frac{\partial \mathbf{g}_{n,2}^*(\beta)}{\partial \beta^T} = \mathcal{E}_n^{[1]}(\beta) + \mathcal{E}_n^{[2]}(\beta) + \mathcal{E}_n^{[3]}(\beta) := \mathcal{E}_n(\beta)$$

where

$$\mathbf{H}_n(\beta) = \sum_{i=1}^n \mathbf{X}_i^T \mathbf{A}_i(\beta)^{1/2} \mathcal{R}_{i-1}^*(\beta)^{-1} \mathbf{A}_i(\beta)^{1/2} \mathbf{X}_i$$

$$\mathbf{B}_n^{[1]}(\beta) = \sum_{i=1}^n \mathbf{X}_i^T \text{diag}\{\mathcal{R}_{i-1}^*(\beta)^{-1} \mathbf{A}_i(\beta)^{-1/2}[\mu_i - \mu_i(\beta)]\} \mathbf{G}_i^{[1]}(\beta) \mathbf{X}_i$$

$$\mathbf{B}_n^{[2]}(\beta) = \sum_{i=1}^n \mathbf{X}_i^T \mathbf{A}_i(\beta)^{1/2} \mathcal{R}_{i-1}^*(\beta)^{-1} \text{diag}\{\mu_i - \mu_i(\beta)\} \mathbf{G}_i^{[2]}(\beta) \mathbf{X}_i$$

$$\mathbf{b}_{n,l}^{[3]}(\beta) = \sum_{i=1}^n \mathbf{X}_i^T \mathbf{A}_i(\beta)^{1/2} \left[\frac{\partial}{\partial \beta_l^T} \mathcal{R}_{i-1}^*(\beta)^{-1}\right] \mathbf{A}_i(\beta)^{-1/2}[\mu_i - \mu_i(\beta)]$$

$$\mathcal{E}_n^{[1]}(\beta) = \sum_{i=1}^n \mathbf{X}_i^T \text{diag}\{\mathcal{R}_{i-1}^*(\beta)^{-1} \mathbf{A}_i(\beta)^{-1/2} \varepsilon_i\} \mathbf{G}_i^{[1]}(\beta) \mathbf{X}_i$$

$$\mathcal{E}_n^{[2]}(\beta) = \sum_{i=1}^n \mathbf{X}_i^T \mathbf{A}_i(\beta)^{1/2} \mathcal{R}_{i-1}^*(\beta)^{-1} \text{diag}(\varepsilon_i) \mathbf{G}_i^{[2]}(\beta) \mathbf{X}_i$$

$$e_{n,l}^{[3]}(\beta) = \sum_{i=1}^{n} \mathbf{X}_i^T \mathbf{A}_i(\beta)^{1/2} \left[ \frac{\partial}{\partial \beta_l^T} \mathcal{R}_{i-1}^*(\beta)^{-1} \right] \mathbf{A}_i(\beta)^{-1/2} \varepsilon_i.$$

Here, the matrices $\mathbf{G}_i^{[1]}(\beta), \mathbf{G}_i^{[2]}(\beta)$ are the same as in [27], i.e.

$$\mathbf{G}_i^{[1]}(\beta) := \mathrm{diag}\left\{ \frac{\mu''(\mathbf{x}_{i1}^T\beta)}{2\mu'(\mathbf{x}_{i1}^T\beta)^{1/2}}, \ldots, \frac{\mu''(\mathbf{x}_{im}^T\beta)}{2\mu'(\mathbf{x}_{im}^T\beta)^{1/2}} \right\}$$

$$\mathbf{G}_i^{[2]}(\beta) := \mathrm{diag}\left\{ -\frac{\mu''(\mathbf{x}_{i1}^T\beta)}{2\mu'(\mathbf{x}_{i1}^T\beta)^{3/2}}, \ldots, -\frac{\mu''(\mathbf{x}_{im}^T\beta)}{2\mu'(\mathbf{x}_{im}^T\beta)^{3/2}} \right\}$$

We denote by $\mathbf{b}_{n,l}^{[3]}(\beta), e_{n,l}^{[3]}(\beta)$ the $l$-th column vectors of $\mathbf{B}_n^{[3]}(\beta)$, respectively $\mathcal{E}_n^{[3]}(\beta)$, for any $1 \leq l \leq p$.

## B  Proof of Lemma 4.6

Clearly, $\eta_n(r) \leq \psi_n(r)$, where

$$\psi_n(r) := \sup_{\beta, \beta' \in B_r} \max_{i \leq n, j \leq m} \left| \frac{\mu'(\mathbf{x}_{ij}^T\beta')}{\mu'(\mathbf{x}_{ij}^T\beta)} - 1 \right|.$$

Let $\beta, \beta' \in B_r$ be arbitrary. By Taylor's formula, there exists $\bar{\beta}_{ij}$ between $\beta$ and $\beta'$ such that $\mu'(\mathbf{x}_{ij}^T\beta') - \mu'(\mathbf{x}_{ij}^T\beta) = \mu''(\mathbf{x}_{ij}^T\bar{\beta}_{ij})\mathbf{x}_{ij}^T(\beta' - \beta)$. By $(AH)$,

$$\left| \frac{\mu'(\mathbf{x}_{ij}^T\beta')}{\mu'(\mathbf{x}_{ij}^T\beta)} - 1 \right| = \left| \frac{\mu''(\mathbf{x}_{ij}^T\bar{\beta}_{ij})}{\mu'(\mathbf{x}_{ij}^T\beta)} \right| |\mathbf{x}_{ij}^T(\beta' - \beta)| \leq C \left| \frac{\mu'(\mathbf{x}_{ij}^T\bar{\beta}_{ij})}{\mu'(\mathbf{x}_{ij}^T\beta)} \right| |\mathbf{x}_{ij}^T(\beta' - \beta)|.$$

Since $|\mathbf{x}_{ij}^T(\beta' - \beta)|^2 \leq \|\mathbf{x}_{ij}^T(\mathbf{H}_n^{\mathrm{indep}})^{-1/2}\|^2 \cdot \|(\mathbf{H}_n^{\mathrm{indep}})^{1/2}(\beta' - \beta)\|^2 \leq \gamma_n^{(0),\mathrm{indep}} \cdot \lambda_{\max}(\mathbf{H}_n^{\mathrm{indep}})r^2 = a_n r^2$, it follows that $\psi_n(r) \leq C\psi_n(r)a_n^{1/2}r + Ca_n^{1/2}r$, i.e.

$$\psi_n(r)(1 - Ca_n^{1/2}r) \leq Ca_n^{1/2}r.$$

By assumption $(K)$, there exist $r_1 > 0$ and $n_1 \geq 1$ such that $a_n^{1/2}r \leq 1/(2C)$ for all $r \in (0, r_1), n \geq n_1$. Hence, $\psi_n(r) \leq Ca_n^{1/2}r$ for all $r \in (0, r_1), n \geq n_1$. □

## C  Proof of Lemma 4.9

We write $\mathbf{H}_n(\beta) - \mathbf{H}_n = \mathbf{H}_n^{[1]}(\beta) + \mathbf{H}_n^{[2]}(\beta) + \mathbf{H}_n^{[3]}(\beta)$, where:

$$\mathbf{H}_n^{[1]}(\beta) = \sum_{i=1}^{n} \mathbf{X}_i^T [\mathbf{A}_i(\beta)^{1/2} - \mathbf{A}_i^{1/2}] \mathcal{R}_{i-1}^*(\beta)^{-1} \mathbf{A}_i(\beta)^{1/2} \mathbf{X}_i$$

$$\mathbf{H}_n^{[2]}(\beta) = \sum_{i=1}^{n} \mathbf{X}_i^T \mathbf{A}_i^{1/2} [\mathcal{R}_{i-1}^*(\beta)^{-1} - (\mathcal{R}_{i-1}^*)^{-1}] \mathbf{A}_i(\beta)^{1/2} \mathbf{X}_i$$

$$\mathbf{H}_n^{[3]}(\beta) = \sum_{i=1}^{n} \mathbf{X}_i^T \mathbf{A}_i^{1/2} (\mathcal{R}_{i-1}^*)^{-1} [\mathbf{A}_i(\beta)^{1/2} - \mathbf{A}_i^{1/2}] \mathbf{X}_i.$$

Let $\lambda$ be an arbitrary $p \times 1$ vector with $\|\lambda\| = 1$. By the Cauchy-Schwartz inequality, $|\lambda^T \mathbf{H}_n^{[1]}(\beta)\lambda| \leq T_1(\beta, \lambda)^{1/2} T_2(\beta, \lambda)^{1/2}$, where

$$T_1(\beta, \lambda) := \sum_{i=1}^n \lambda^T \mathbf{X}_i^T \mathbf{A}_i(\beta)^{1/2} \mathcal{R}_{i-1}^*(\beta)^{-1} \mathbf{A}_i(\beta)^{1/2} \mathbf{X}_i \lambda \leq \max_{i \leq n} \lambda_{\max}[\mathcal{R}_{i-1}^*(\beta)^{-1}]$$

$$\max_{i \leq n} \lambda_{\max}^2[\mathbf{A}_i^{-1/2} \mathbf{A}_i(\beta)^{1/2}] \cdot \lambda^T \mathbf{H}_n^{\mathrm{indep}} \lambda \leq C \pi_n(r) \lambda_{\max}(\mathbf{H}_n^{\mathrm{indep}})$$

$$T_2(\beta, \lambda) := \sum_{i=1}^n \lambda^T \mathbf{X}_i^T [\mathbf{A}_i(\beta)^{1/2} - \mathbf{A}_i^{1/2}] \mathcal{R}_{i-1}^*(\beta)^{-1} [\mathbf{A}_i(\beta)^{1/2} - \mathbf{A}_i^{1/2}] \mathbf{X}_i \lambda$$

$$\leq \max_{i \leq n} \lambda_{\max}[\mathcal{R}_{i-1}^*(\beta)^{-1}] \max_{i \leq n} \lambda_{\max}^2[\mathbf{A}_i^{-1/2} \mathbf{A}_i(\beta)^{1/2} - \mathbf{I}]$$

$$\max_{i \leq n} \lambda_{\max}^2[\mathbf{A}_i^{-1/2} \mathbf{A}_i(\beta)^{1/2}] \cdot \lambda^T \mathbf{H}_n^{\mathrm{indep}} \lambda \leq C \pi_n(r) \eta_n^2(r) \lambda_{\max}(\mathbf{H}_n^{\mathrm{indep}}).$$

For estimating the terms above, we used $(R')$, (48) and Lemma 4.6. Hence,

$$|\lambda^T \mathbf{H}_n^{[1]}(\beta)\lambda| \leq C \lambda_{\max}(\mathbf{H}_n^{\mathrm{indep}}) \pi_n(r) \eta_n(r).$$

Note that $|\lambda^T \mathbf{H}_n^{[2]}(\beta)\lambda| \leq T_1'(\lambda)^{1/2} T_2'(\beta, \lambda)^{1/2}$ where

$$T_1'(\lambda) := \sum_{i=1}^n \lambda^T \mathbf{X}_i^T \mathbf{A}_i^{1/2} (\mathcal{R}_{i-1}^*)^{-1} \mathbf{A}_i^{1/2} \mathbf{X}_i \lambda \leq C \lambda_{\max}(\mathbf{H}_n^{\mathrm{indep}})$$

$$T_2'(\beta, \lambda) := \sum_{i=1}^n \lambda^T \mathbf{X}_i^T \mathbf{A}_i(\beta)^{1/2} (\mathcal{R}_{i-1}^*)^{-1/2} [(\mathcal{R}_{i-1}^*)^{1/2} \mathcal{R}_{i-1}^*(\beta)^{-1} (\mathcal{R}_{i-1}^*)^{1/2} - \mathbf{I}]^2$$

$$(\mathcal{R}_{i-1}^*)^{-1/2} \mathbf{A}_i(\beta)^{1/2} \mathbf{X}_i \lambda$$

$$\leq \max_{i \leq n} \lambda_{\max}[(\mathcal{R}_{i-1}^*)^{1/2} \mathcal{R}_{i-1}^*(\beta)^{-1} (\mathcal{R}_{i-1}^*)^{1/2} - \mathbf{I}]^2 \max_{i \leq n} \lambda_{\max}[(\mathcal{R}_{i-1}^*)^{-1}]$$

$$\max_{i \leq n} \lambda_{\max}[\mathbf{A}_i^{-1/2} \mathbf{A}_i(\beta)^{1/2}] \cdot \lambda^T \mathbf{H}_n^{\mathrm{indep}} \lambda$$

$$\leq C \rho_n^2(r) \lambda_{\max}(\mathbf{H}_n^{\mathrm{indep}}).$$

Hence

$$|\lambda^T \mathbf{H}_n^{[2]}(\beta)\lambda| \leq \lambda_{\max}(\mathbf{H}_n^{\mathrm{indep}}) \rho_n(r).$$

Similarly, one ca prove that:

$$|\lambda^T \mathbf{H}_n^{[3]}(\beta)\lambda| \leq \lambda_{\max}(\mathbf{H}_n^{\mathrm{indep}}) \eta_n(r).$$

Since $\pi_n(r) \geq 1$, by $(C_1')$ and $(C_2)$,

$$\lim_{r \to 0} \limsup_{n \to \infty} \alpha_n^{-1/2-\delta} \sup_{\beta \in B_r} \sup_{\|\lambda\|=1} |\lambda^T \mathbf{H}_n^{[k]}(\beta)\lambda| = 0 \quad \text{a.s,}$$

for $k = 1, 2, 3$. It follows that

$$\lim_{r \to 0} \limsup_{n \to \infty} \alpha_n^{-1/2-\delta} \sup_{\beta \in B_r} \sup_{\|\lambda\|=1} |\lambda^T [\mathbf{H}_n(\beta) - \mathbf{H}_n]\lambda| = 0 \quad \text{a.s.}$$

□

# D  Proof of Lemma 4.10

Let $\lambda$ be an arbitrary $p \times 1$ vector with $\|\lambda\| = 1$. Using the fact that $\mathrm{diag}(\mathbf{v})\mathbf{D}\mathbf{u} = \mathbf{D}\mathrm{diag}(\mathbf{u})\mathbf{v}$, for any vectors $\mathbf{u}, \mathbf{v}$ and for any diagonal matrix $\mathbf{D}$, we write:

$$\lambda^T \mathbf{B}_n^{[1]}(\beta)\lambda = \sum_{i=1}^n \lambda^T \mathbf{X}_i^T \mathbf{G}_i^{[1]}(\beta)\mathrm{diag}(\mathbf{X}_i\lambda)\mathcal{R}_{i-1}^*(\beta)^{-1}\mathbf{A}_i(\beta)^{-1/2}[\mu_i - \mu_i(\beta)]$$

$$\lambda^T \mathbf{B}_n^{[2]}(\beta)\lambda = \sum_{i=1}^n \lambda^T \mathbf{X}_i^T \mathbf{A}_i(\beta)^{1/2}\mathcal{R}_{i-1}^*(\beta)^{-1}\mathrm{diag}(\mathbf{X}_i\lambda)\mathbf{G}_i^{[2]}(\beta)[\mu_i - \mu_i(\beta)].$$

By the Cauchy-Schwartz inequality, $|\lambda^T \mathbf{B}_n^{[1]}(\beta)\lambda| \leq I_1(\beta,\lambda)^{1/2} I_2(\beta)^{1/2}$, where

$$\begin{aligned}
I_1(\beta,\lambda) &:= \sum_{i=1}^n \lambda^T \mathbf{X}_i \mathbf{G}_i^{[1]}(\beta)\mathrm{diag}(\mathbf{X}_i\lambda)\mathcal{R}_{i-1}^*(\beta)^{-1}\mathrm{diag}(\mathbf{X}_i\lambda)\mathbf{G}_i^{[1]}(\beta)\mathbf{X}_i\lambda \\
&\leq \max_{i\leq n}\lambda_{\max}[\mathcal{R}_{i-1}^*(\beta)^{-1}]\max_{i\leq n}\lambda_{\max}[\mathrm{diag}^2(\mathbf{X}_i\lambda)]\max_{i\leq n}\lambda_{\max}^2[\mathbf{A}_i^{-1/2}\mathbf{G}_i^{[1]}(\beta)] \\
&\quad \lambda^T \mathbf{H}_n^{\mathrm{indep}}\lambda \\
&\leq C\pi_n(r)a_n\lambda_{\max}(\mathbf{H}_n^{\mathrm{indep}}) \\
I_2(\beta) &:= \sum_{i=1}^n [\mu_i - \mu_i(\beta)]^T \mathbf{A}_i(\beta)^{-1/2}\mathcal{R}_{i-1}^*(\beta)^{-1}\mathbf{A}_i(\beta)^{-1/2}[\mu_i - \mu_i(\beta)] \\
&\leq \max_{i\leq n}\lambda_{\max}[\mathcal{R}_{i-1}^*(\beta)^{-1}]\max_{i\leq n}\lambda_{\max}^2[\mathbf{A}_i(\bar{\beta}_i)^{1/2}\mathbf{A}_i(\beta)^{-1}\mathbf{A}_i(\bar{\beta}_i)^{1/2}] \\
&\quad \max_{i\leq n}\lambda_{\max}[\mathbf{A}_i^{-1}\mathbf{A}_i(\bar{\beta}_i)] \cdot (\beta - \beta_0)^T \mathbf{H}_n^{\mathrm{indep}}(\beta - \beta_0) \\
&\leq C\pi_n(r)(\beta-\beta_0)^T \mathbf{H}_n^{\mathrm{indep}}(\beta-\beta_0) \\
&\leq Cr^2\pi_n(r)\lambda_{\max}(\mathbf{H}_n^{\mathrm{indep}}). \quad (53)
\end{aligned}$$

For estimating the term $I_2(\beta)$, we used Lemma 4.6, (48), and the Taylor's formula: $\mu_i(\beta) - \mu_i = \mathbf{A}_i(\bar{\beta}_i)\mathbf{X}_i(\beta - \beta_0)$, where $\bar{\beta}_i$ is between $\beta$ and $\beta_0$. For $I_1(\beta, \lambda)$, we used (48), $(AH)$, Lemma 4.6, and the fact that

$$\lambda_{\max}[\mathrm{diag}^2(\mathbf{X}_i\lambda)] \leq a_n \quad \text{for all } i \leq n, \qquad (54)$$

(To prove (54), note that $|\mathbf{x}_{ij}^T\lambda|^2 \leq \lambda_{\max}(\mathbf{H}_n^{\mathrm{indep}})|\mathbf{x}_{ij}^T(\mathbf{H}_n^{\mathrm{indep}})^{-1}\mathbf{x}_{ij}| \leq \lambda_{\max}(\mathbf{H}_n^{\mathrm{indep}})\gamma_n^{(0),\mathrm{indep}} = a_n$.) Therefore

$$|\lambda^T \mathbf{B}_n^{[1]}(\beta)\lambda| \leq Cra_n^{1/2}\pi_n(r)\lambda_{\max}(\mathbf{H}_n^{\mathrm{indep}}).$$

Using condition $(C_1)$, it follows that

$$\lim_{r\to 0}\limsup_{n\to\infty}\alpha_n^{-1/2-\delta}\sup_{\beta\in B_r}\sup_{\|\lambda\|=1}|\lambda^T \mathbf{B}_n^{[1]}(\beta)\lambda| = 0 \quad \text{a.s.}$$

The term $\mathbf{B}_n^{[2]}(\beta)$ is treated by similar methods. More precisely, $|\lambda^T \mathbf{B}_n^{[2]}(\beta)\lambda| \leq I_1'(\beta,\lambda)^{1/2}I_2'(\beta)^{1/2}$, where

$$I_1'(\beta,\lambda) := \sum_{i=1}^n \lambda^T \mathbf{X}_i^T \mathbf{A}_i(\beta)^{1/2}\mathcal{R}_{i-1}^*(\beta)^{-1}\mathrm{diag}(\mathbf{X}_i\lambda)\mathbf{G}_i^{[2]}(\beta)\mathbf{A}_i(\beta)^{1/2}\mathcal{R}_{i-1}^*$$

$$\mathbf{A}_i(\beta)^{1/2}\mathbf{G}_i^{[2]}(\beta)\mathrm{diag}(\mathbf{X}_i\lambda)\mathcal{R}_{i-1}^*(\beta)^{-1}\mathbf{A}_i(\beta)^{1/2}\mathbf{X}_i\lambda$$
$$\leq C\pi_n^2(r)\max_{i\leq n}\lambda_{\max}^2[\mathbf{G}_i^{[2]}(\beta)\mathbf{A}_i(\beta)\mathbf{G}_i^{[2]}(\beta)]\max_{i\leq n}\lambda_{\max}[\mathrm{diag}^2(\mathbf{X}_i\lambda)]$$
$$\max_{i\leq n}\lambda_{\max}[\mathbf{A}_i^{-1}\mathbf{A}_i(\beta)]\cdot\lambda^T\mathbf{H}_n^{\mathrm{indep}}\lambda$$
$$\leq C\pi_n^2(r)a_n\lambda_{\max}(\mathbf{H}_n^{\mathrm{indep}})$$
$$I_2'(\beta) := \sum_{i=1}^n[\mu_i-\mu_i(\beta)]^T\mathbf{A}_i(\beta)^{-1/2}(\mathcal{R}_{i-1}^*)^{-1}\mathbf{A}_i(\beta)^{-1/2}[\mu_i-\mu_i(\beta)]$$
$$\leq Cr^2\lambda_{\max}(\mathbf{H}_n^{\mathrm{indep}}).$$

To estimate $I_1'(\beta,\lambda)$ we used (48), (54), Lemma 4.6, $(R')$ and $(AH)$. The estimate of $I_2'(\beta)$ was obtained similarly to (53), using $(R')$.

By condition $(C_1)$, we conclude that:
$$\lim_{r\to 0}\limsup_{n\to\infty}\alpha_n^{-1/2-\delta}\sup_{\beta\in B_r}\sup_{\|\lambda\|=1}|\lambda^T\mathbf{B}_n^{[2]}(\beta)\lambda| = 0 \quad \text{a.s.}$$

To treat the term which involves $\mathbf{B}_n^{[3]}(\beta)$, we note that $|\lambda^T\mathbf{B}_n^{[3]}(\beta)\lambda|^2 \leq \sum_{l=1}^p |\lambda^T\mathbf{b}_{n,l}^{[3]}|^2$, where $\mathbf{b}_{n,l}^{[3]}(\beta)$ denotes the $l$-th column of $\mathbf{B}_n^{[3]}(\beta)$. By Theorem 9.2, [22], we have:
$$\frac{\partial}{\partial\beta_l}\mathcal{R}_{i-1}^*(\beta)^{-1} = -\mathcal{R}_{i-1}^*(\beta)^{-1}\left[\frac{\partial}{\partial\beta_l}\mathcal{R}_{i-1}^*(\beta)\right]\mathcal{R}_{i-1}^*(\beta)^{-1}. \qquad (55)$$

Now, for any $l\in\{1,\ldots,p\}$ fixed, we have: $|\lambda^T\mathbf{b}_{n,l}^{[3]}| \leq I_{1,l}''(\beta,\lambda)^{1/2}I_2(\beta)^{1/2}$, where $I_2(\beta)$ is as above and
$$I_{1,l}''(\beta,\lambda) := \sum_{i=1}^n \lambda^T\mathbf{X}_i^T\mathbf{A}_i(\beta)^{1/2}\mathcal{R}_{i-1}^*(\beta)^{-1}\left[\frac{\partial}{\partial\beta_l}\mathcal{R}_{i-1}^*(\beta)\right]\mathcal{R}_{i-1}^*(\beta)^{-1}\left[\frac{\partial}{\partial\beta_l}\mathcal{R}_{i-1}^*(\beta)\right]$$
$$\mathcal{R}_{i-1}^*(\beta)^{-1}\mathbf{A}_i(\beta)^{1/2}\mathbf{X}_i\lambda$$
$$\leq \max_{i\leq n}\lambda_{\max}^2\left[\frac{\partial}{\partial\beta_l}\mathcal{R}_{i-1}^*(\beta)\right]\max_{i\leq n}\lambda_{\max}^3[\mathcal{R}_{i-1}^*(\beta)^{-1}]\max_{i\leq n}\lambda_{\max}[\mathbf{A}_i^{-1/2}\mathbf{A}_i(\beta)^{1/2}]$$
$$\lambda^T\mathbf{H}_n^{\mathrm{indep}}\lambda$$
$$\leq C\pi_n^3(r)q_n^2(r)\lambda_{\max}(\mathbf{H}_n^{\mathrm{indep}}),$$

using (48) and Lemma 4.6. Hence,
$$|\lambda^T\mathbf{b}_{n,l}^{[3]}| \leq Cr\pi_n^2(r)q_n(r)\lambda_{\max}(\mathbf{H}_n^{\mathrm{indep}}).$$

Finally, $(C_3)$ implies that $\lim_{r\to 0}\limsup_{n\to\infty}|\lambda^T\mathbf{b}_{n,l}^{[3]}| = 0$ a.s., and therefore
$$\lim_{r\to 0}\limsup_{n\to\infty}\alpha_n^{-1/2-\delta}\sup_{\beta\in B_r}\sup_{\|\lambda\|=1}|\lambda^T\mathbf{B}_n^{[3]}(\beta)\lambda| = 0 \quad \text{a.s.}$$

□

# E   Proof of Lemma 4.11

Let $\lambda$ be an arbitrary $p \times 1$ vector with $\|\lambda\| = 1$. We first treat the terms $\mathcal{E}_n^{[1]}(\beta)$ and $\mathcal{E}_n^{[2]}(\beta)$. Using the fact that $\operatorname{diag}(\mathbf{v})\mathbf{D}\mathbf{u} = \mathbf{D}\operatorname{diag}(\mathbf{u})\mathbf{v}$, for any vectors $\mathbf{u}, \mathbf{v}$ and for any diagonal matrix $\mathbf{D}$, we write $\lambda^T \mathcal{E}_n^{[1]}(\beta)\lambda = \sum_{k=1}^{3} U_n^{[k]}(\beta, \lambda)$ and $\lambda^T \mathcal{E}_n^{[2]}(\beta)\lambda = \sum_{k=4}^{6} \mathbf{U}_n^{[k]}(\beta, \lambda)$, where

$$U_n^{[1]}(\beta,\lambda) = \sum_{i=1}^{n} \lambda^T \mathbf{X}_i^T \mathbf{G}_i^{[1]} \operatorname{diag}(\mathbf{X}_i\lambda) \mathcal{R}_{i-1}^*(\beta)^{-1} \mathbf{A}_i^{-1/2} \varepsilon_i$$

$$U_n^{[2]}(\beta,\lambda) = \sum_{i=1}^{n} \lambda^T \mathbf{X}_i^T \mathbf{G}_i^{[1]}(\beta) \operatorname{diag}(\mathbf{X}_i\lambda) \mathcal{R}_{i-1}^*(\beta)^{-1} (\mathbf{A}_i^{-1/2}(\beta) - \mathbf{A}_i^{-1/2}) \varepsilon_i$$

$$U_n^{[3]}(\beta,\lambda) = \sum_{i=1}^{n} \lambda^T \mathbf{X}_i^T (\mathbf{G}_i^{[1]}(\beta) - \mathbf{G}_i^{[1]}) \operatorname{diag}(\mathbf{X}_i\lambda) \mathcal{R}_{i-1}^*(\beta)^{-1} \mathbf{A}_i^{-1/2} \varepsilon_i$$

$$U_n^{[4]}(\beta,\lambda) = \sum_{i=1}^{n} \lambda^T \mathbf{X}_i^T \mathbf{A}_i^{1/2} \mathcal{R}_{i-1}^*(\beta)^{-1} \operatorname{diag}(\mathbf{X}_i\lambda) \mathbf{G}_i^{[2]} \varepsilon_i$$

$$U_n^{[5]}(\beta,\lambda) = \sum_{i=1}^{n} \lambda^T \mathbf{X}_i^T (\mathbf{A}_i^{1/2}(\beta) - \mathbf{A}_i^{1/2}) \mathcal{R}_{i-1}^*(\beta)^{-1} \operatorname{diag}(\mathbf{X}_i\lambda) \mathbf{G}_i^{[2]}(\beta) \varepsilon_i$$

$$U_n^{[6]}(\beta,\lambda) = \sum_{i=1}^{n} \lambda^T \mathbf{X}_i^T \mathbf{A}_i^{1/2} \mathcal{R}_{i-1}^*(\beta)^{-1} \operatorname{diag}(\mathbf{X}_i\lambda) (\mathbf{G}_i^{[2]}(\beta) - \mathbf{G}_i^{[2]}) \varepsilon_i$$

Note that $\{U_n^{[k]}(\beta, \lambda), \mathcal{F}_n\}_n$ is a martingale (with respect to $P_{\beta_0}$), and hence $\{\sup_{\beta \in B_r} \sup_{\|\lambda\|=1} |U_n^{[k]}(\beta, \lambda)|, \mathcal{F}_n\}_n$ is a submartingale (with respect to $P_{\beta_0}$), for any $r > 0$.

In what follows, we will prove that for any $r > 0$, there exists a constant $C > 0$ (depending on $r$) such that

$$E \sup_{\beta \in B_r} \sup_{\|\lambda\|=1} |U_n^{[k]}(\beta,\lambda)| \leq C, \quad \forall n \geq 1. \tag{56}$$

By the martingale convergence theorem (Theorem 2.5, [10]), it will follow that $\{\sup_{\beta \in B_r} \sup_{\|\lambda\|=1} |U_n^{[k]}(\beta, \lambda)|\}_n$ converges a.s. Using the fact that $\alpha_n \to \infty$, we obtain that: for any $r > 0$,

$$\lim_{n \to \infty} \alpha_n^{-1/2-\delta} \sup_{\beta \in B_r} \sup_{\|\lambda\|=1} |U_n^{[k]}(\beta,\lambda)| = 0 \quad \text{a.s.,} \quad k=1,\ldots,6,$$

from which it will follow that

$$\lim_{n \to \infty} \alpha_n^{-1/2-\delta} \sup_{\beta \in B_r} \sup_{\|\lambda\|=1} |\lambda^T \mathcal{E}_n^{[k]}(\beta)\lambda| = 0 \quad \text{a.s.,} \quad k=1,2.$$

We now turn to the proof of (56). By the Cauchy-Schwartz inequality $|U_n^{[1]}(\beta,\lambda)| \leq J_1(\beta,\lambda)^{1/2} J_2^{1/2}$, where $J_2 := \sum_{i=1}^{n} \varepsilon_i^T \mathbf{A}_i^{-1} \varepsilon_i$ and

$$J_1(\beta,\lambda) := \sum_{i=1}^{n} \lambda^T \mathbf{X}_i^T \mathbf{G}_i^{[1]} \operatorname{diag}(\mathbf{X}_i\lambda) \mathcal{R}_{i-1}^*(\beta)^{-2} \operatorname{diag}(\mathbf{X}_i\lambda) \mathbf{G}_i^{[1]} \mathbf{X}_i\lambda$$

$$
\begin{aligned}
&\leq \max_{i\leq n}\lambda_{\max}^2[\mathcal{R}_{i-1}^*(\beta)^{-1}]\max_{i\leq n}\{\lambda_{\max}[\mathrm{diag}^2(\mathbf{X}_i\lambda)]\}\max_{i\leq n}\lambda_{\max}^2[\mathbf{A}_i^{-1/2}\mathbf{G}_i^{[1]}]\\
&\quad\lambda^T\mathbf{H}_n^{\mathrm{indep}}\lambda\\
&\leq C\pi_n^2(r)a_n\lambda_{\max}(\mathbf{H}_n^{\mathrm{indep}}).
\end{aligned}
$$

For the estimation of the term $J_1(\beta,\lambda)$, we used (48), (54) and $(AH)$.

Note that $E(\varepsilon_i^T\mathbf{A}_i^{-1}\varepsilon_i)=\mathrm{tr}\overline{\mathbf{R}}_i^{(c)}=m$ for all $i$, and hence

$$E(J_2) = E\left(\sum_{i=1}^n \varepsilon_i^T\mathbf{A}_i^{-1}\varepsilon_i\right) = mn. \tag{57}$$

We conclude that

$$
\begin{aligned}
E\sup_{\beta\in B_r}\sup_{\|\lambda\|=1}|U_n^{[1]}(\beta,\lambda)| &\leq \{E\sup_{\beta\in B_r}\sup_{\|\lambda\|=1}J_1(\beta,\lambda)\}^{1/2}\{E(J_2)\}^{1/2}\\
&\leq C\{E[\pi_n^2(r)]na_n\lambda_{\max}(\mathbf{H}_n^{\mathrm{indep}})\}^{1/2}.
\end{aligned}
$$

Similarly, we find the upper bound $C\{E[\pi_n^2(r)]na_n^2\lambda_{\max}(\mathbf{H}_n^{\mathrm{indep}})\}^{1/2}$ for the term $E\sup_{\beta\in B_r}\sup_{\|\lambda\|=1}|U_n^{[k]}(\beta,\lambda)|$, with $k=2,3$. For this, we use the following fact: for any $r>0$

$$\sup_{\beta\in B_r}\max_{i\leq n,j\leq m}\left\{\left|\frac{\mu''(\mathbf{x}_{ij}^T\beta)}{\mu'(\mathbf{x}_{ij}^T\beta)^{1/2}}-\frac{\mu''(\mathbf{x}_{ij}^T\beta_0)}{\mu'(\mathbf{x}_{ij}^T\beta_0)^{1/2}}\right|\mu'(\mathbf{x}_{ij}^T\beta_0)^{-1/2}\right\}\leq Cra_n^{1/2}.$$

(This inequality can be proved using Taylor's formula.)

Relation (56) with $k=1,2,3$, follows from $(C_4)$.

We now treat the term involving $U_n^{[4]}(\beta,\lambda)$. By the Cauchy-Schwartz inequality, $|U_n^{[4]}(\beta,\lambda)|\leq J_3(\beta,\lambda)^{1/2}J_4(\beta,\lambda)^{1/2}$, where

$$
\begin{aligned}
J_3(\beta,\lambda) &:= \sum_{i=1}^n \lambda^T\mathbf{X}_i\mathbf{A}_i^{1/2}\mathcal{R}_{i-1}^*(\beta)^{-1}\mathrm{diag}^2(\mathbf{X}_i\lambda)\mathcal{R}_{i-1}^*(\beta)^{-1}\mathbf{A}_i^{1/2}\mathbf{X}_i\lambda\\
&\leq \max_{i\leq n}\{\lambda_{\max}[\mathrm{diag}^2(\mathbf{X}_i\lambda)]\}\max_{i\leq n}\lambda_{\max}[\mathcal{R}_{i-1}^*(\beta)^{-2}]\cdot\lambda^T\mathbf{H}_n^{\mathrm{indep}}\lambda\\
&\leq C\pi_n^2(r)a_n\lambda_{\max}(\mathbf{H}_n^{\mathrm{indep}})\\
J_4(\beta,\lambda) &:= \sum_{i=1}^n \varepsilon_i^T\left(\mathbf{G}_i^{[2]}\right)^2\varepsilon_i \leq \max_{i\leq n}[\mathbf{A}_i^{1/2}\mathbf{G}_i^{[2]}]\sum_{i=1}^n\varepsilon_i^T\mathbf{A}_i^{-1}\varepsilon_i\\
&\leq C\sum_{i=1}^n\varepsilon_i^T\mathbf{A}_i^{-1}\varepsilon_i,
\end{aligned}
$$

and we used (48) and $(AH)$. By (57), $E\sup_{\beta\in B_r}\sup_{\|\lambda\|=1}J_4(\beta,\lambda)\leq Cn$, and hence

$$
\begin{aligned}
E\sup_{\beta\in B_r}\sup_{\|\lambda\|=1}|U_n^{[4]}(\beta,\lambda)| &\leq \{E\sup_{\beta\in B_r}\sup_{\|\lambda\|=1}J_3(\beta,\lambda)\}^{1/2}\{E\sup_{\beta\in B_r}\sup_{\|\lambda\|=1}J_4(\beta,\lambda)\}^{1/2}\\
&\leq C\{E[\pi_n^2(r)]na_n\lambda_{\max}(\mathbf{H}_n^{\mathrm{indep}})\}^{1/2}.
\end{aligned}
$$

Similarly, we find the upper bound $C\{E[\pi_n^2(r)]na_n^2\lambda_{\max}[\mathbf{H}_n^{\text{indep}}]\}^{1/2}$ for the term $E\sup_{\beta\in B_r}\sup_{\|\lambda\|=1}|U_n^{[k]}(\beta,\lambda)|$, with $k=5,6$, using the fact that: $\forall r>0$

$$\sup_{\beta\in B_r}\max_{i\leq n, j\leq m}\left\{\left|\frac{\mu''(\mathbf{x}_{ij}^T\beta)}{\mu'(\mathbf{x}_{ij}^T\beta)^{3/2}}-\frac{\mu''(\mathbf{x}_{ij}^T\beta_0)}{\mu'(\mathbf{x}_{ij}^T\beta_0)^{3/2}}\right|\mu'(\mathbf{x}_{ij}^T\beta_0)^{-1/2}\right\}\leq Cra_n^{1/2}.$$

Relation (56) with $k=4,5,6$, follows from $(C_4)$.

It remains to treat the term $\mathcal{E}_n^{[3]}(\beta)$. Note that $|\lambda^T\mathcal{E}_n^{[3]}\lambda|^2\leq\sum_{l=1}^{p}|\lambda^T\mathbf{e}_{n,l}^{[3]}(\beta)|^2$, where $\mathbf{e}_{n,l}^{[3]}(\beta)$ denotes the $l$-th column of $\mathcal{E}_n^{[3]}(\beta)$. For any $1\leq l\leq p$, we define $U_{n,l}^{[7]}(\beta,\lambda)=\lambda^T\mathbf{e}_{n,l}^{[3]}(\beta)$.

Note that $\{U_{n,l}^{[7]}(\beta,\lambda),\mathcal{F}_n\}_n$ is a martingale (with respect to $P_{\beta_0}$), and hence $\{\sup_{\beta\in B_r}\sup_{\|\lambda\|=1}|U_{n,l}^{[7]}(\beta,\lambda)|,\mathcal{F}_n\}_n$ is a submartingale, for any $r>0$. Using the same argument as above, to conclude that for any $r>0$,

$$\lim_{n\to\infty}\alpha_n^{-1/2-\delta}\sup_{\beta\in B_r}\sup_{\|\lambda\|=1}|\lambda^T\mathcal{E}_n^{[3]}(\beta)\lambda|=0\quad\text{a.s.,}$$

it suffices to show that

$$\lim_{n\to\infty}\alpha_n^{-1/2-\delta}\sup_{\beta\in B_r}\sup_{\|\lambda\|=1}|U_{n,l}^{[7]}(\beta,\lambda)|=0\quad\text{a.s..}$$

For this, it is enough to show that: for any $r>0$, there exists a constant $C>0$ such that

$$E\sup_{\beta\in B_r}\sup_{\|\lambda\|=1}|U_{n,l}^{[7]}(\beta,\lambda)|\leq C,\quad\forall n\geq 1. \tag{58}$$

By the Cauchy-Schwartz inequality, $|U_{n,l}^{[7]}(\beta,\lambda)|\leq J_5(\beta,\lambda)^{1/2}J_6(\beta,\lambda)^{1/2}$, where

$$\begin{aligned}
J_5(\beta,\lambda) &:= \sum_{i=1}^{n}\lambda^T\mathbf{X}_i^T\mathbf{A}_i(\beta)^{1/2}\mathcal{R}_{i-1}^*(\beta)^{-1}\left[\frac{\partial}{\partial\beta_l}\mathcal{R}_{i-1}^*(\beta)\right]\mathcal{R}_{i-1}^*(\beta)^{-2}\\
&\qquad\left[\frac{\partial}{\partial\beta_l}\mathcal{R}_{i-1}^*(\beta)\right]\mathcal{R}_{i-1}^*(\beta)^{-1}\mathbf{A}_i(\beta)^{1/2}\mathbf{X}_i\lambda\\
&\leq \max_{i\leq n}\lambda_{\max}[\mathcal{R}_{i-1}^*(\beta)^{-4}]\max_{i\leq n}\lambda_{\max}^2\left[\frac{\partial}{\partial\beta_l}\mathcal{R}_{i-1}^*(\beta)\right]\max_{i\leq n}\lambda_{\max}[\mathbf{A}_i^{-1}\mathbf{A}_i(\beta)]\\
&\qquad \lambda^T\mathbf{H}_n^{\text{indep}}\lambda\\
&\leq C\pi_n^4(r)q_n^2(r)\lambda_{\max}(\mathbf{H}_n^{\text{indep}})\\
J_6(\beta,\lambda) &:= \sum_{i=1}^{n}\varepsilon_i^T\mathbf{A}_i(\beta)^{-1}\varepsilon_i\leq\max_{i\leq n}\lambda_{\max}[\mathbf{A}_i^{1/2}\mathbf{A}_i(\beta)^{-1}\mathbf{A}_i^{1/2}]\sum_{i=1}^{n}\varepsilon_i^T\mathbf{A}_i^{-1}\varepsilon_i\\
&\leq C\sum_{i=1}^{n}\varepsilon_i^T\mathbf{A}_i^{-1}\varepsilon_i,
\end{aligned}$$

and we used (55), $(AH)$ and (48) and Lemma 4.6. By (57), $E\sup_{\beta\in B_r}\sup_{\|\lambda\|=1}J_6(\beta,\lambda)\leq Cn$, and hence

$$E\sup_{\beta\in B_r}\sup_{\|\lambda\|=1}|U^{[7]}_{n,l}(\beta,\lambda)| \leq \{E\sup_{\beta\in B_r}\sup_{\|\lambda\|=1}J_5(\beta,\lambda)\}^{1/2}\{E\sup_{\beta\in B_r}\sup_{\|\lambda\|=1}J_6(\beta,\lambda)\}^{1/2}$$
$$\leq C\{E[\pi_n^4(r)q_n^2(r)]n\lambda_{\max}(\mathbf{H}_n^{\text{indep}})\}^{1/2}.$$

Relation (58) follows by condition $(C_5)$. $\square$